\documentclass[10pt,fleqn]{article}
\usepackage{mathptmx}

\input{epsf.tex}
\usepackage{latexsym,amsfonts,amssymb,epsfig,verbatim}
\usepackage{amsmath,amsthm,amssymb,latexsym,graphics,textcomp}
\usepackage{eucal,eufrak}
\usepackage{graphicx}
\usepackage{url}

\theoremstyle{plain}
\newtheorem{theorem}{Theorem}[section]
\newtheorem{lemma}[theorem]{Lemma}
\newtheorem{corollary}[theorem]{Corollary}

\newtheorem{proposition}[theorem]{Proposition}

\newtheorem*{theoremdsub}{Theorem \ref{domainsub}}

\theoremstyle{definition}
\newtheorem*{remark}{Remark}

\newtheorem{question}{Question}

\def\1{\mathbf 1}

\def\ext{\mathrm{ext}}

\def\UML{\mathcal{UML}}
\def\EL{\mathcal{EL}}
\def\ML{\mathcal{ML}}
\def\MF{\mathcal{MF}}
\def\PML{\mathbb P \mathcal{ML}}
\def\PMF{\mathbb P\mathcal{MF}}
\def\WH{\mathfrak{H}}
\def\T{\mathrm{Teich}}
\def\Mod{\mathrm{Mod}}
\def\Homeo{\mathrm{Homeo}}

\def\A{\mathfrak A}
\def\C{\mathcal C}
\def\F{\mathcal F}
\def\M{\mathcal M}
\def\N{\mathcal N}
\def\S{\mathcal S}
\def\T{\mathcal T}
\def\Q{\mathcal Q}

\def\L{\mathcal L}

\def\diam{\mathrm{diam}}

\def\d{\mathrm{d}}

\def\Dom{\mathrm{Dom}}

\def\Tel{\mathcal T_{\mathrm{el}}}
\def\Rel{\mathbb R_{\mathrm{el}}}

\def\co{\colon\thinspace}

\def\g{\mathcal G}

\begin{document}

\title{\textbf{Shadows of mapping class groups: \\ capturing convex cocompactness}}
\author{Richard P. Kent IV and  Christopher J. Leininger \thanks{The second author was supported by an N.S.F. postdoctoral fellowship.}}
\maketitle

\textit{for Dick and Geri}


\section{Introduction}\label{introsect}

A Kleinian group $\Gamma$ is a discrete subgroup of $\mathrm{PSL}_2(\mathbb C)$.  When non-elementary, such a group possesses a unique non-empty minimal closed invariant subset $\Lambda_\Gamma$ of the Riemann sphere, called the limit set. A Kleinian group acts properly discontinuously on the complement $\Delta_\Gamma$ of $\Lambda_\Gamma$ and so this set is called the domain of discontinuity.

Such a group is said to be convex cocompact if it acts cocompactly on the convex hull $\mathrm H_\Gamma$ in $\mathbb H^3$ of its limit set $\Lambda_\Gamma$. This is equivalent to the condition that an orbit of $\Gamma$ is quasi-convex in $\mathbb H^3$---or that the orbit defines a quasi-isometric embedding $\Gamma \to \mathbb H^3$. Equivalent to each of these is the property that every limit point of $\Gamma$ is conical, and still another definition is that $\Gamma$ has a compact Kleinian manifold---meaning that $\Gamma$ acts cocompactly on $\mathbb H^3 \cup \Delta_\Gamma$. We refer the reader to \cite{bowditchgeom} and the references therein for the history of these notions and the proof of their equivalence (see also \cite{swenson}). 

Let $S$ denote an oriented complete hyperbolic surface of finite area, $\Mod(S) = \pi_0(\Homeo^+(S))$ its group of orientation preserving self--homeomorphisms up to isotopy, and $\T(S)$ the Teichm\"uller space of $S$ equipped with Teichm\"uller's metric.

The mapping class group $\Mod(S)$ acts on Teichm\"uller space $\T(S)$ by isometries, and W. Thurston discovered a $\Mod(S)$--equivariant compactification of $\T(S)$ by an ideal sphere, the sphere of compactly supported projective measured laminations $\PML(S)$.
J. McCarthy and A. Papadopoulos have shown that a subgroup $G$ of $\Mod(S)$ has a well defined limit set $\Lambda_G$, although it need not be unique or minimal, and that there is a certain enlargement $Z\Lambda_G$ of $\Lambda_G$ on whose complement $G$ acts properly discontinuously \cite{mccarthypapa}.  So such a group has a domain of discontinuity $\Delta_G= \PML(S) - Z\Lambda_G$.

In general, the limit set of a subgroup of $\Mod(S)$ has no convex hull to speak of, as there are pairs in $\PML(S)$ that are joined by no geodesic in $\T(S)$. Nevertheless, if every pair of points in $\Lambda_G$ are the negative and positive directions of a geodesic in $\T(S)$, one can define the weak hull $\WH_G$ of $\Lambda_G$ to be the union of all such geodesics. This is precisely what B. Farb and L. Mosher do in \cite{FMcc}, where they develop a notion of convex cocompact mapping class groups.
They prove the following

\begin{theorem} {\bf (Farb--Mosher)} \label{famoconv}
Given a finitely generated subgroup $G$ of $\Mod(S)$, the following statements are equivalent:
\begin{itemize}
\item Some orbit of $G$ is quasi-convex in $\T(S)$.
\item Every orbit of $G$ is quasi-convex in $\T(S)$.
\item $G$ is 
hyperbolic and there is a $G$--equivariant embedding $\partial f\co \partial G \to \PML(S)$ with image $\Lambda_G$ such that the weak hull $\WH_G$ of $\Lambda_G$ is defined; the action of $G$ on $\WH_G$ is cocompact; and, if $f\co G \rightarrow \WH_G$ is any $G$--equivariant map, then $f$ is a quasi-isometry and the following map is continuous:
\[
\overline{f} = f \cup \partial f\co G \cup \partial G \to \T(S) \cup \PML(S).
\]
\end{itemize}
\end{theorem}

\noindent A finitely generated subgroup of $\Mod(S)$ is said to be \textbf{convex cocompact} if it satisfies one of these conditions.\\

\noindent The interest in convex cocompact Kleinian groups is due in part to combined work of A. Marden \cite{Mar} and D. Sullivan \cite{Sul} that implies that such groups are precisely those that are quasiconformally stable, meaning that small perturbations of the identity representation are induced by quasiconformal conjugacies.

The allure of convex cocompact mapping class groups is of a manifestly different nature:  Farb and Mosher have shown that when $S$ is closed, convex cocompactness for a subgroup $G < \Mod(S)$ is implied by the $\delta$--hyperbolicity of the associated $\pi_1(S)$--extension of $G$---see \cite{FMcc}.
Moreover, in very recent work, U. Hamenst\"adt has shown that these are equivalent \cite{hamenstadt}.
In particular, if there is a finite $\mathrm{K}(G,1)$, and an embedding $G \to \Mod(S)$ whose image is not convex cocompact and yet whose
non-identity elements are all pseudo-Anosov, then the associated surface group extension is a group with no Baumslag--Solitar subgroups, a finite
Eilenberg--Mac Lane space, and which fails to be hyperbolic.  This would provide a counterexample to a question of M. Gromov---see \cite{klsurvey} and \cite{FMcc}.

For more on the geometry of these extensions and related groups, we refer the reader to \cite{FMcc}, \cite{FMII}, \cite{Mitra1}, \cite{Mitra2}, \cite{Mitra3}, \cite{hypbyhyp}, and \cite{combinationtheorem}.\\

\noindent Our purpose here is to strengthen the analogy between convex cocompact Kleinian groups and their cousins in the mapping class group. Our first main result is the following

\begin{theorem}\label{main} Given a finitely generated subgroup $G$ of $\Mod(S)$, the following statements are equivalent:
\begin{itemize}
\item $G$ is convex cocompact.
\item The weak hull $\WH_G$ is defined and $G$ acts cocompactly on $\WH_G$.
\item Every limit point of $G$ is conical.
\item $G$ acts cocompactly on $\T(S) \cup \Delta_G$.
\end{itemize}
\end{theorem}

\begin{remark} The definition of $\WH_G$ used here is more general than that described above and is defined for any infinite irreducible $G$, see Section \ref{trianglesect}.
\end{remark}

That $G$ need only act cocompactly on $\WH_G$ to be convex cocompact follows quickly from the fact, proven in Section \ref{trianglesect}, that geodesic triangles lying in a thick part of $\T(S)$ are thin in the sense of $\delta$--hyperbolic metric spaces: if $G$ acts cocompactly on $\WH_G$ it is coarsely dense therein and the weak hull lies in a thick part of $\T(S)$; the thin triangle condition on $\WH_G$ implies that it is quasi-convex [Theorem \ref{convexhull}], and a $G$--orbit is quasi-convex as a result.  That triangles lying in a thick part are thin relies on H. Masur's Asymptotic Rays Theorem \cite{masuruniquely} and Y. Minsky's Contraction Theorem \cite{minskycrelle}.

\begin{figure}[h!]
\begin{center}
\input{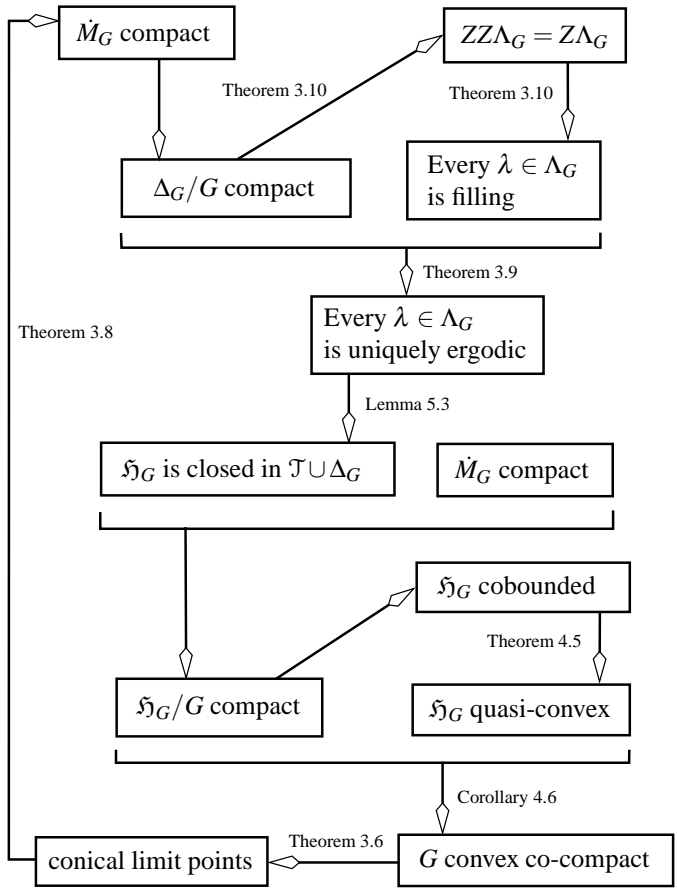}
\end{center}
\caption{Convex cocompactness, Kleinian manifolds, and conical limit points.}
\label{logictwo}
\end{figure}

If $\Gamma$ is a Kleinian group, a limit point for $\Gamma$ is said to be \textbf{conical} if every geodesic in $\mathbb H^3$ terminating there has a neighborhood that intersects a $\Gamma$--orbit in an infinite set. In the mapping class group, the definition of conical requires some care as there are points in the boundary of Thurston's compactification of $\T(S)$ that are not limits of Teichm\"uller geodesic rays.  Nonetheless, points exhibiting this behavior are irrelevant by Masur's Two Boundaries Theorem \cite{twoboundaries}, and it is easily seen that convex cocompact groups have all limit points conical [Theorem \ref{ccconical}].

With the aid of F. Bonahon's work on geodesic currents \cite{bonahon}, the arguments given by McCarthy and Papadopoulos to prove that $G$ acts properly discontinuously on $\Delta_G$ can be extended to prove proper discontinuity on $\T(S) \cup \Delta_G$.
We write $\dot{M}_G = \linebreak (\T(S) \cup \Delta_G) /G$ and refer to this as the \textbf{Kleinian manifold} for $G$.
Along with certain length and intersection number comparisons along Teichm\"uller geodesic rays, these extended arguments prove that if all limit points are conical, then $\dot{M}_G$ is compact [Theorem \ref{conicalcompact}]. The only remaining implication is that having a compact Kleinian manifold implies convex cocompactness.

Minsky's Bounded Geometry Theorem \cite{minskybound} for a doubly degenerate hyperbolic $3$--manifold with a type preserving homeomorphism to $S \times \mathbb R$ says that the length of the shortest geodesic of such a manifold is bounded below if and only if the Masur--Minsky subsurface projection coefficients of its ending laminations are uniformly bounded above.  K. Rafi has proven the analog of this theorem for Teichm\"uller geodesics \cite{rafi}: namely, a geodesic lies in a thick part of $\T(S)$ if and only if all of the subsurface projection coefficients of its defining laminations are uniformly bounded.

The set $Z\Lambda_G$ is the set of laminations having zero intersection number with some lamination in $\Lambda_G$. The set $Z Z \Lambda_G$ is the set of laminations having zero intersection number with some element of $Z\Lambda_G$. We may continue this procedure to obtain a sequence of sets $Z^{(n)}\Lambda_G$. When a subgroup $G$ of $\Mod(S)$ acts cocompactly on $\Delta_G$, $Z\Lambda_G$ is stable under this operation of taking zero loci and every lamination in $\Lambda_G$ is filling [Theorem \ref{zz}].  
A cocompact action on $\Delta_G$, in conjunction with Rafi's bounded geometry theorem for Teichm\"uller geodesics, implies that, in fact, every lamination in $\Lambda_G$ is uniquely ergodic [Theorem \ref{domainsub}]. Such groups always have weak hulls that are closed in $\T(S) \cup \Delta_G$ [Lemma \ref{thehullisclosed}] and compactness of $\WH_G/G$ follows from compactness of $\dot M_G$. The logic of the proof of Theorem \ref{main} is depicted in Figure \ref{logictwo}.

Theorem \ref{domainsub} provides much stronger information than is needed to prove Theorem \ref{main}.  We state it here as it may be of independent interest.
\begin{theoremdsub} Let $G$ be a subgroup of $\Mod(S)$. If $\Delta_G \neq \emptyset$ and $G$ acts cocompactly on $\Delta_G$, then every lamination in $\Lambda_G$ is uniquely ergodic, $Z\Lambda_G = \Lambda_G$, and $\WH_G$ is defined and cobounded. Furthermore, $G$ has a finite index subgroup all of whose non-identity elements are pseudo-Anosov.
\end{theoremdsub}

\noindent An earlier proof that convex cocompact mapping class groups have compact Kleinian manifolds mirrored the proof in the Kleinian group setting and revealed that weak hulls lying in a thick part of Teichm\"uller space have closest points projections with contraction properties similar to convex hulls in $\mathbb H^3$, generalizing the quasi-projection theorems of Minsky---we have preserved this projection theorem in Section \ref{hullprojsect}.\\

\noindent An obstacle to shining light on $\Mod(S)$ presents itself when one has taken a point of view based on the analogy between $\T(S)$ and $\mathbb H^3$: the Teichm\"uller space with the Teichm\"uller metric is not hyperbolic in any reasonable sense of the word \cite{masurclass,masurwolf} (nor is it hyperbolic with any reasonable $\Mod(S)$--invariant metric \cite{brockfarb}).  
Indeed, even if the map sending a subgroup $G$ of $\Mod(S)$ to its orbit in $\T(S)$ is a quasi-isometric embedding, $G$ need not be convex cocompact---not even when $G$ is cyclic \cite{MMunstable}. On the other hand, $\Mod(S)$ acts by isometries on W. Harvey's complex of curves $\C(S)$,  which is $\delta$--hyperbolic by a celebrated theorem of H. Masur and Y. Minsky \cite{MM1,bowditch}. Illuminating $\Mod(S)$ from this vantage point has some advantages over the view from $\T(S)$---as well as disadvantages due to the fact that $\C(S)$ fails to be locally compact.  Our second main theorem is the following.

\begin{theorem}\label{complextheorem} A finitely generated subgroup $G$ of $\Mod(S)$ is convex cocompact if and only if sending $G$ to an orbit in the complex of curves defines a quasi-isometric embedding $G \to \C(S)$.
\end{theorem}

\begin{remark}
This theorem was independently discovered by U. Hamenst\"adt \cite{hamenstadt}.
\end{remark}

\begin{figure}
\begin{center}
\input{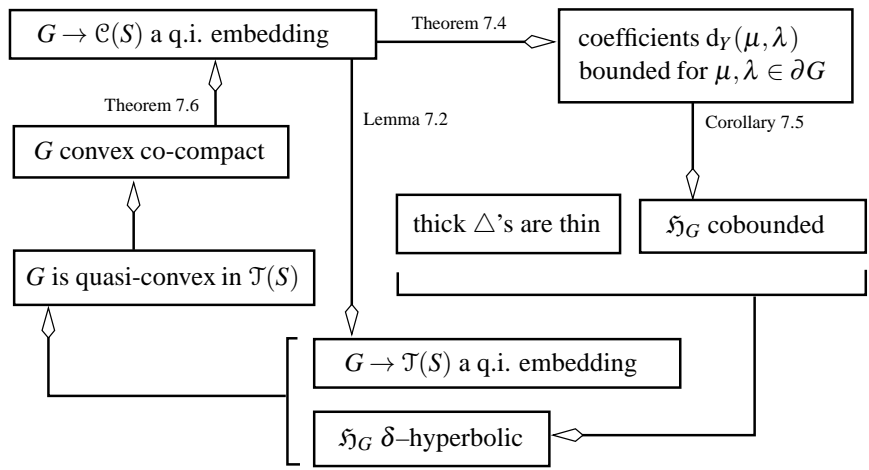}
\end{center}
\caption{Quasi-isometrically embedding in $\C(S)$ and convex cocompactness.}
\label{logicthree}
\end{figure}

The proof that convex cocompact groups have this property is a novel application of Minsky's Contraction Theorem combined with the fact that $\C(S)$ is quasi-isometric to the electric Teichm\"uller space $\Tel(S)$.

Given a quasi-isometric embedding $G\to \C(S)$, we obtain from general principles that $G$ is $\delta$--hyperbolic and that the given map extends continuously to a map 
\[
G\cup \partial G \to \C(S) \cup \partial \C(S)
\]
whose restriction to $\partial G$ is an embedding. The boundary of the complex of curves is naturally parameterized by the space $\EL(S)$ of potential ending laminations for geometrically infinite hyperbolic $3$--manifolds homeomorphic to $S \times \mathbb R$, by a theorem of E. Klarreich \cite{klarreich}. The space $\EL(S)$ sits naturally---as the set of filling laminations---in the quotient of $\PML(S)$ obtained by forgetting transverse measures.

Using hyperbolicity of $\C(S)$ and Masur and Minsky's Bounded Geodesic Image Theorem \cite{MM2}, we are able to uniformly bound the projection coefficients for endpoints in $\partial \C(S)$ of bi-infinite geodesics in $G$. Rafi's Bounded Geometry Theorem and a proposition of Klarreich allow us to lift $\partial G$ to the limit set $\Lambda_G$  and demonstrate that the weak hull $\WH_G$ is defined and cobounded. The fact that triangles in a thick part are thin again tells us that $\WH_G$ is $\delta$--hyperbolic. The quasi-isometric embedding $G \to \C(S)$ yields a quasi-isometric embedding $G \to \WH_G$ and we conclude that $G$ is quasi-convex in $\T(S)$ by hyperbolicity of the hull. See Figure \ref{logicthree}.\\

\noindent In \cite{FMcc}, Farb and Mosher prove that when $S$ is closed, a free subgroup of $\Mod(S)$ is convex cocompact if and only if the associated surface-by-free group is $\delta$--hyperbolic.  Such subgroups are called \textbf{Schottky}. In \cite{hypbyhyp}, Mosher proves that for any finite set of independent pseudo-Anosov mapping classes $\varphi_1,\, \ldots \, ,\varphi_n$---meaning that the fixed points are pairwise distinct---there is a natural number $m$ such that the surface-by-free group associated to $\varphi_1^m,\, \ldots \, , \varphi_n^m$ is $\delta$--hyperbolic.  This demonstrates the abundance of Schottky groups. Theorem \ref{complextheorem} yields a new proof of this fact (without the closed hypothesis on $S$).

\begin{theorem}[Abundance of Schottky groups \cite{FMcc}]\label{abundance} Given a finite set of independent pseudo-Anosov mapping classes $\varphi_1, \, \ldots \, ,\varphi_n$, there is a number $\ell$ so that for all natural numbers $m > \ell$, the group generated by $\varphi_1^m, \, \ldots \, ,\varphi_n^m$ is Schottky.
\end{theorem}
\begin{proof} We refer the reader to Section \ref{backgroundsect} for terminology.

Let $\mathfrak{G}_1,\, \ldots \, ,\mathfrak{G}_n$ denote the Cayley graphs of $\langle \varphi_1 \rangle,\, \ldots \, ,\langle \varphi_n \rangle$, respectively.
Fix $\alpha \in \C(S)$ and $\langle \varphi_i \rangle$--equivariant embeddings $\mathfrak{G}_i \to \C(S)$ by sending each vertex to the associated point of the orbit of $\alpha$ and sending edges to geodesics joining the images of their endpoints.
By Theorem 4.6 of \cite{MM1}, these are all $(K_0,C_0)$--quasi-isometric embeddings for some $K_0 \geq 1$ and $C_0 \geq 0$.
Moreover, since $\varphi_1, \, \ldots \, , \varphi_n$ are independent, all the endpoints of these quasi-geodesic rays in $\partial \C(S)$ are distinct.


For any $m > 0$, we consider the abstract free group $\langle \varphi_1^m,\, \ldots \, ,\varphi_n^m \, | \, \mbox{---} \, \rangle$ equipped with the metric $m\, \d$, where $\d$ is the word metric with respect to the $\varphi_i^m$.
Note that we naturally obtain a metric on the associated Cayley graph $\mathfrak{G}(m)$.
There is a canonical isometric embedding $\mathfrak{G}_i \to \mathfrak{G}(m)$ restricting to the ``identity'' on $\langle \varphi_i^m \rangle$.
Our chosen embeddings of $\mathfrak{G}_i \to \C(S)$ induces a $\langle \varphi_1^m,\, \ldots \, ,\varphi_n^m \, | \, \mbox{---} \, \rangle$--equivariant map $\mathfrak{G}(m) \to \C(S)$.

Now, given two quasi-geodesic rays $\tau_1$ and $\tau_2$ with a common origin in a $\delta$--hyperbolic metric space $\mathcal X$ representing distinct points in $\partial \mathcal X$, the union $\tau_1 \cup \tau_2$ (suitably parameterized) is a quasi-geodesic line with constants depending only on the $\tau_i$, their quasi-geodesic constants, and $\delta$.

Fix $K \geq K_0 \geq 1$ and $C \geq C_0 \geq 1$ quasi-geodesic constants for the embeddings of $\mathfrak{G}_i$ into $\C(S)$ and all quasi-geodesic lines obtained as unions of distinct rays in unions $\mathfrak{G}_i \cup \mathfrak{G}_j$ (via the embeddings into $\C(S)$).

Since $\C(S)$ is $\delta$--hyperbolic for some $\delta$, there is an $R \geq 1$ such that any $(K,C,R)$--local-quasi-geodesic  is a quasi-geodesic, see Th\'eor\`eme 1.4 of \cite{coor}.
If $m \geq R$, then every geodesic segment of length $R$ through $\1$ in $\mathfrak{G}(m)$ is contained in a union $\mathfrak{G}_i \cup \mathfrak{G}_j$.
In fact, such a segment is contained in a union of two geodesic rays from $\1$ contained in $\mathfrak{G}_i$ and $\mathfrak{G}_j$ respectively and by the choice of $K$ and $C$, this segment is sent to a $(K,C)$--quasi-geodesic segment in $\C(S)$.

As any geodesic segment in $\mathfrak{G}(m)$ may be translated to a segment through $\1$, we conclude that every geodesic in $\mathfrak{G}(m)$ is sent to a $(K,C,R)$--local-quasi-geodesic in $\C(S)$, and that $\langle \varphi_1^m,\, \ldots \, ,\varphi_n^m \, | \, \mbox{---} \, \rangle$ is quasi-isometrically embedded by its orbit in $\C(S)$ provided $m \geq R$.
In particular, $\langle \varphi_1^m, \, \ldots \, , \varphi_n^m \rangle \to \Mod(S)$ is injective, and has convex cocompact image by Theorem \ref{complextheorem}.
\end{proof}

\begin{remark} Although Farb and Mosher work with closed surfaces in considering convex cocompactness, their definitions carry over verbatim for the case of finite area hyperbolic surfaces.  Theorem \ref{famoconv} also easily extends to this setting.  A quick way to see this (that requires no verification of Farb-Mosher's techniques) is to observe that the Teichm\"uller spaces of punctured surfaces isometrically embed in the Teichm\"uller spaces of closed surfaces with nice mapping class group equivariance properties (via appropriate branched covers).

In addition, using this method, Farb and Mosher's proof of corollary \ref{abundance} easily implies the corollary for finite area hyperbolic surfaces.  Furthermore, we note that this trick allows any example of a convex cocompact subgroups of finite area hyperbolic surfaces to be promoted virtually to examples in closed surface mapping class groups.
\end{remark}

\noindent \textbf{Acknowledgments.} The authors thank Jeff Brock, Jason DeBlois, Moon Duchin, Benson Farb, Cameron Gordon, Yair Minsky, Lee Mosher, Alan Reid, and Peter Storm for useful conversations.  We also thank the referee for several helpful comments and suggestions.

\section{Background} \label{backgroundsect}

\subsection{Coarse geometry}\label{coarsesect}

Given metric spaces $\mathcal X$ and $\mathcal Y$ and constants $K \geq 1$ and $C \geq 0$, a map $f\co \mathcal X \to \mathcal Y$ is a \textbf{$(K,C)$--quasi-isometric embedding} if
\[
\frac{1}{K}\d_{\mathcal X}(a,b) - C \leq \d_{\mathcal Y}(f(a),f(b)) \leq K \d_{\mathcal X}(a,b) + C
\]
for all $a$ and $b$ in $\mathcal X$, and a \textbf{$(K,C)$--quasi-isometry} if its image is $A$--dense for some $A$. Such a map is said to be \textbf{$K$--bi-Lipschitz} if $C=0$.

A map $f \co \mathcal X \to \mathcal Y$ is said to be \textbf{$(K,C)$--coarsely--Lipschitz} if
\[
\d_{\mathcal Y}(f(a),f(b)) \leq K \d_{\mathcal X}(a,b) + C
\]
for all $a$ and $b$ in $\mathcal X$, and \textbf{$K$--Lipschitz} if $(K,0)$--coarsely--Lipschitz.

A map from an interval in $\mathbb R$ or $\mathbb Z$ to a metric space $\mathcal X$ is a \textbf{$(K,C)$--quasi-geodesic} if it is a $(K,C)$--quasi-isometric embedding and a \textbf{$(K,C,R)$--local--quasi-geodesic} if its restriction to any interval of length $R$ is a $(K,C)$--quasi-isometric embedding.

If $\mathcal X$ is a geodesic metric space, $\mathcal Y \subset \mathcal X$, then $\mathcal Y$ is said to be \textbf{$A$--quasi-convex} if every geodesic joining two points in $\mathcal Y$ is contained in the $A$--neighborhood $\N_A(\mathcal Y)$ of $\mathcal Y$.

Given a finitely generated group $G$ with finite generating set $\mathcal U$, let $\d_{\mathcal U}$ denote the induced word metric.
A \textbf{geodesic} $\g$ in $(G,\d_{\mathcal U})$ is a $(1,0)$--quasi-geodesic defined on an interval $\mathbb{I} \subset \mathbb{Z}$.
We represent $\g$ by a sequence of group elements $\g = \{ h_j \}_{j \in \mathbb{I}}$ and emphasize that the defining characteristic of being a $(1,0)$--quasi-geodesic is that
\[ \d_{\mathcal U}(h_j,h_k) = |j-k|. \]

Note that, given any point $h_i \in \g$, we can translate $\g$ to a geodesic through the identity $\1$ taking $h_i$ to $\1$, namely
\[
h_i^{-1}(\g) = \{ h_i^{-1} h_j \}_{j \in \mathbb I}\, .
\]

\subsection{$\delta$--hyperbolic spaces}\label{hypsect}

We refer the reader to \cite{coor} and \cite{BH} for more on hyperbolic metric spaces and coarse geometry.

A geodesic triangle in a metric space is \textbf{$\delta$--thin} if each of its sides is contained in the $\delta$--neighborhood of the union of the other two sides.

A geodesic metric space is \textbf{$\delta$--hyperbolic} in the sense of M. Gromov and J. Cannon if every geodesic triangle is $\delta$--thin.

Let $\mathcal X$ be a metric space. Given $x$, $y$, and $z$ in $\mathcal X$, the \textbf{Gromov product of $y$ and $z$} with respect to $x$ is defined to be
\[
(y \cdot z)_x = \frac 1 2 (\d_{\mathcal X}(y,x) + \d_{\mathcal X}(z,x) - \d_{\mathcal X}(y,z)).
\]
Fix a basepoint $x$ in $\mathcal X$. A sequence $\{x_n\}$ \textbf{converges at infinity} in $\mathcal X$ if 
\[
\lim_{m,n \to \infty} (x_m \cdot x_n)_x = \infty
\]
and two sequences $\{x_n\}$ and $\{y_m\}$ are \textbf{equivalent} if 
\[
\lim_{m,n \to \infty} (y_m \cdot x_n)_x = \infty.
\]
If $\mathcal X$ is a $\delta$--hyperbolic geodesic metric space, the \textbf{Gromov boundary $\partial \mathcal X$} of $\mathcal X$ is the set of equivalence classes of sequences in $\mathcal X$ that converge at infinity.  The set $\overline{\mathcal X} = \mathcal X \cup \partial \mathcal X$ admits a natural topology in which a sequence $\{x_n\}$ in $\mathcal X$ converges to a point $y = \{y_n\}$ in $\partial \mathcal X$ if and only if $\{x_n\}$ is equivalent to $y$. 

A geodesic ray based at $x$ uniquely determines a point in $\partial \mathcal X$ given by any sequence of points on the ray that converges at infinity.

\subsection{Teichm\"uller theory}\label{teichsect}

We refer the reader to \cite{ahlforsqc,abikoff,gardiner} for more on quasiconformal mappings and Teichm\"uller theory.  In what follows, unless otherwise stated, all Riemann surfaces and \linebreak complex/conformal structures are of finite analytic type, and hyperbolic surfaces and structures are complete with finite area.

Let $X$ be a Riemann surface homeomorphic to $S$.
We view $X$ as either a complex $1$--manifold, or an oriented hyperbolic $2$--manifold---the Uniformization Theorem permits us to change this view at will.
A \textbf{marking} of $X$ is an orientation preserving homeomorphism $f \co S \to X$, and the Teichm\"uller space $\T(S)$ of $S$ is the set of equivalence classes of marked Riemann surfaces $f \co S \to X$.
The equivalence relation is defined by declaring $f_1 \co S \to X_1$ to be equivalent to $f_2 \co S \to X_2$ if $f_2 \circ f_1^{-1}$ is isotopic to an isomorphism of Riemann surfaces.

Abusing notation, we often refer to a Riemann surface $X$ as a point in Teichm\"uller space, by which we mean the equivalence class of $X$ implicitly marked by some homeomorphism.

We may also think of $\T(S)$ as the space of complex, conformal, or hyperbolic structures on $S$, up to isotopy, as such a structure is induced on $S$ by pulling back via the marking.

Let $X_1$ and $X_2$ be two Riemann surfaces equipped with markings $f_1\co S \to X_1$ and $f_2\co S \to X_2$.  A homeomorphism $f\co X_1 \to X_2$ is \textbf{$K$--quasiconformal} if it is absolutely continuous on lines and $|f_{\bar z}| \leq k |f_z|$ in every local coordinate $z$ where $k= \linebreak (K-1)/(K+1) < 1$. The minimal value of $K$ for which $f$ is $K$--quasiconformal is the \textbf{dilatation} of $f$ and is denoted $K[f]$.  The \textbf{Teichm\"uller distance} between $X_1$ and $X_2$ is defined to be
\[
\d_\T(X_1,X_2) = \frac 1 2 \inf \log K[f]
\]
where the infimum is taken over all quasiconformal maps $f$ isotopic to $f_2 \circ f_1^{-1}$. There is a unique extremal quasiconformal map $X_1 \to X_2$  realizing the above distance, called the \textbf{Teichm\"uller mapping}.

A \textbf{holomorphic quadratic differential} $q$ on $X$ is an assignment of a holomorphic function $\varphi(z)$ to each local coordinate $z$ such that for two coordinates $z_1$ and $z_2$, 
\[
\varphi_1(z_1) (\d z_1/\d z_2)^2 = \varphi_2 (z_2).
\]
We say that $q$ is \textbf{integrable} if $\int_X|q| < \infty$.
We equip the vector space of integrable holomorphic quadratic differentials $\Q(X)$ with the norm $\|\ \ \|=\int_X|\ \ |$.
Varying $X$ over $\T(S)$ and assembling the vector spaces $\Q(X)$ one obtains a vector bundle $\Q(S)$ over $\T(S)$.
We let $\Q^1(S)$ denote the associated unit sphere bundle over $\T(S)$ and $\Q^*(S)$ the complement of the zero section of $\Q(S) \to \T(S)$.
Throughout the remainder of this paper, we will refer to a point of $\Q^*(S)$ simply as a quadratic differential, with the holomorphic, integrability, and non-zero conditions implicit.  We denote a quadratic differential by $q$ or sometimes $(X,q)$ if we wish to emphasize the underlying complex structure $X \in \T(S)$.

A \textbf{natural parameter} at a point $p$ associated to a quadratic differential $q = \varphi(z) \, \d z^2$ is given by
\[
\zeta(w) = \int_{z_0}^w \sqrt{\varphi(z)} \ \d z \, ,
\] 
where $z_0 = z(p)$.  Away from the zeroes of $q$, $\zeta$ is a \textit{bona fide} coordinate for which $q$ takes the particularly simple form $q = \d \zeta^2$.

Pulling back the horizontal and vertical foliations of $\mathbb C$ via a  natural parameter $\zeta$ yields a pair of transverse measured (singular) foliations $\F_-(q)$ and $\F_+(q)$ associated to $q$, where the transverse measures are the pullbacks of the vertical and horizontal total variations in $\mathbb C$, respectively. A branched flat metric, called the \textbf{$q$--metric}, is obtained from $\zeta$ by pulling back the Euclidean metric from $\mathbb C$.  The total area of $X$ with respect to the $q$-metric is $\|q\|$.

The Teichm\"uller mapping between Riemann surfaces may be described explicitly in terms of a holomorphic quadratic differential.
Fix such a differential $q$ with natural parameter $\zeta$ and a number $0 \leq k < 1$.
The \textbf{$(q,k)$--Teichm\"uller deformation} of $X$ is the Riemann surface $X'$ with underlying surface $X$ and complex structure given by the local parameter $\zeta' = (\zeta + k\overline{\zeta})/(1-k)$.
As $X'$ is determined by $X$, $q$, and $k$, we let $(X,q,k)$ denote $X'$.
The quadratic differential $q' = (d\zeta')^2/\|(d \zeta')^2\|$ is called the \textbf{terminal differential} of the deformation.
The horizontal and vertical measured foliations of the terminal differential are given by $\mathcal{F}_-(q') = K^{-1/2}\mathcal{F}_-(q)$ and $\mathcal{F}_+(q') = K^{1/2}\mathcal{F}_+(q)$, where $K = (1+k)(1-k)$.
The ``identity'' $X \to X'$ is the Teichm\"uller mapping in its homotopy class, and $\d_{\mathcal{T}}(S,S') = \frac{1}{2} \log K$.
Teichm\"uller's Theorem asserts that for any $X \in \T(S)$ the map from the unit ball $B_1 \subset \Q(X)$ to $\mathcal{T}(S)$ given by $q \mapsto (X,q,\|q\|)$ is a homeomorphism.

Given a quadratic differential $q$ on $X$, define the \textbf{Teichm\"uller geodesic} $\tau = \linebreak \tau_q \co \mathbb{R} \to \mathcal{T}(S)$ as follows.
For $t \in \mathbb{R}$, let $s_t$ be given by $t = \frac{1}{2} \log((1 + s_t)/(1 - s_t))$, and define
\[ \tau(t) = \left\{ \begin{array}{ll}
(X,q,s_t) & \mbox{for } t > 0\\
(X,-q,-s_t) & \mbox{otherwise}\\ \end{array} \right. \]
The parameter is chosen so that $\tau$ is a geodesic parameterized by arc--length with respect to the Teichm\"uller metric---note that the terminal differential at time $t$ has horizontal and vertical measured foliations $\mathcal{F}_-(q_t) = e^{-t}\mathcal{F}_-(q)$ and $\mathcal{F}_+(q_t) = e^t\mathcal{F}_+(q)$.
We refer to the projective classes of the horizontal and vertical foliations of $q$ as the horizontal and vertical foliations of $\tau$.

The mapping class group $\Mod(S)$ acts on $\T(S)$ by pulling back hyperbolic structures.  This is a properly discontinuous action by isometries of the Teichm\"uller metric and the quotient $\M(S)$ is the \textbf{moduli space} of $S$.

The \textbf{$\epsilon$--thick part} of $\T(S)$ is the set of hyperbolic 
structures on $S$ for which the length of the shortest geodesic is greater than $\epsilon > 0$.  A 
set is said to be \textbf{$\epsilon$--cobounded} if it lies in the 
$\epsilon$--thick part of $\T(S)$ and \textbf{cobounded} if it is $\epsilon$--cobounded for some $\epsilon > 0$. D. Mumford's Compactness Criterion \cite{mumford} says that a set is $\epsilon$--cobounded for some $\epsilon$ if and only if it projects to a precompact set in $\M(S)$, see also \cite{notesonnotes}.

\subsection{Laminations and foliations}\label{lamssect}

For a nice survey of geodesic laminations, see \cite{bonahonsurvey}.

Fix a hyperbolic metric on $S$.  A \textbf{geodesic lamination} on $S$ is a nonempty compact subset of $S$ that is a union of pairwise disjoint simple complete local geodesics on $S$, called the \textbf{leaves} of the lamination.  By a simple complete local geodesic, we mean the image of an injective locally isometric immersion of $\mathbb R$. A \textbf{transverse measure} on a geodesic lamination is an assignment of a Radon measure to each arc (and so each curve) $\alpha$ transverse to the lamination such that (1) the measure on a subarc $\alpha'$ of an arc $\alpha$ is the restriction to $\alpha'$ of the measure on $\alpha$ and (2) so that two arcs are assigned the same measure if they are isotopic through arcs transverse to the lamination.  A \textbf{measured lamination} $\lambda$ is a geodesic lamination $|\lambda|$ called the \textbf{support} of $\lambda$, along with a transverse measure.  We will further always require that our transverse measures have \textbf{full support}:  if the intersection of a transverse arc $\alpha$ with $|\lambda|$ is non-empty, then the measure on $\alpha$ is non-zero (equivalently, the support of the measure on $\alpha$ is exactly $\alpha \cap |\lambda|$).

The set of measured laminations admits a natural topology, see \cite{bonahonsurvey}, and the resulting space is denoted $\ML(S)$. The set of non-zero laminations up to scale, the \textbf{projective measured laminations}, is denoted $\PML(S)$.  We let $\UML(S)$ denote the quotient of the set of measured laminations obtained by forgetting transverse measures. The space $\ML(S)$ depends on the choice of hyperbolic metric, but there is a canonical identification between spaces obtained by different metrics.  

Certain definitions for (or properties of) $\lambda \in \ML(S)$ may depend only on its projective class $[\lambda]$ or its support $|\lambda|$.  In these cases, we will freely apply such definitions (or refer to such properties) for $[\lambda]$ or $|\lambda|$, as is appropriate.
Occasionally, we confuse a measured lamination $\lambda$ with its projective class $[\lambda]$, and even its support $|\lambda|$, referring to all three simply as $\lambda$ when the meaning is clear.

Let $\S$ be the set of isotopy classes of essential simple closed curves on $S$. Essential here means homotopically nontrivial and nonperipheral.  Taking geodesic representatives with transverse counting measures, we identify $\mathcal{S}$ with a subset of $\ML(S)$. The geometric intersection number $i(\, \cdot \, ,\, \cdot \,)\co \S \times \S \to \mathbb R_+=[0,\infty)$ extends naturally to a continuous function
\[
i(\, \cdot \, ,\, \cdot \,)\co \ML(S) \times \ML(S) \to \mathbb R_+ \, .
\]
When $\lambda$ is a measured lamination and $\alpha$ is a simple closed curve, $i(\lambda,\alpha)$ is simply the total mass of transverse measure on $\alpha$ assigned by $\lambda$.  The injection $\S \to \ML(S)$ remains injective upon descending to the quotients $\S \to \PML(S)$ and $\S \to \UML(S)$.

A measured lamination $\lambda$ is said to be \textbf{filling} if it intersects non-trivially any measured lamination whose support is different from that of $\lambda$.  This property for $\lambda$ depends only on $|\lambda|$.

There is a closely related theory of (singular) \textbf{measured foliations} on $S$.  There is a space $\MF(S)$ of (measure classes of) measured foliations, an intersection function $i:\MF(S) \times \MF(S) \to {\mathbb R}$, and a space of such foliations up to scale, $\PMF(S)$. The horizontal and vertical foliations of a holomorphic quadratic differential $q \in \Q(S)$ are examples of transversely measured foliations.  In fact, it is a theorem of J. Hubbard and H. Masur \cite{hubbardmasur} and S. Kerckhoff \cite{kerckhoff}  that for any fixed complex structure, each (measure class) of measured foliation appears as the vertical foliation of a unique holomorphic quadratic differential.

As discussed in the previous subsection, for any $X \in \T(S)$, Teichm\"uller's Theorem provides a homeomorphism from the open unit ball in $\Q(X)$ to $\T(S)$.  Therefore, the closed unit ball serves as a visual compactification of $\mathcal{T}(S)$. Identifying a quadratic differential with the projective class of its vertical foliation, we interpret this as a compactification by $\PMF(S)$.
This is \textbf{Teichm\"uller's compactification} of $\mathcal{T}(S)$.

We have chosen to work primarily with measured laminations rather than foliations, and so refer the reader to \cite{FLP} for a detailed treatment of the latter. However, we need to know that there is a ``dictionary'' between $\ML(S)$ and $\MF(S)$---see G. Levitt's paper \cite{levitt} for details. Given a measured foliation $F$, there is an associated lamination $\lambda_F$, obtained from $F$ by straightening the leaves of $F$.  The foliation $F$ can be recovered from $\lambda_F$ by a certain collapsing procedure applied to the complementary regions of $\lambda_F$ (when $\lambda_F$ has isolated leaves, one must first replace such a leaf by an annulus foliated by parallel copies of the core).  We call $\lambda_F$ the \textbf{lamination associated to} $F$, and $F$ the \textbf{foliation associated to} $\lambda_F$.

The identification $\ML(S) \cong \MF(S)$ is completely natural.  Specifically, we have $i(F_1,F_2) = i(\lambda_{F_1},\lambda_{F_2})$ and for every $t > 0$ and $F \in \MF(S)$, one has $\lambda_{tF} = t \lambda_F$.  Moreover, for any $g \in \Mod(S)$, we have $g \lambda_F = \lambda_{gF}$.

Given a Teichm\"uller geodesic $\tau$ defined by a holomorphic quadratic differential $q$ with horizontal and vertical foliations $\F_-$ and $\F_+$, we call the associated laminations $\lambda_- = \lambda_{\F_-}$ and $\lambda_+ = \lambda_{\F_+}$ (or their projective classes) the \textbf{negative and positive directions} of $\tau$, respectively.  If $\tau$ is a geodesic ray, the lamination $\lambda_+$ associated to $\F_+$ is called the \textbf{direction} of $\tau$.  The boundary of Teichm\"uller's compactification of $\T(S)$ is thus identified with $\PML(S)$ as the directions of rays emanating from a point $X \in \T(S)$.

If two projective measured geodesic laminations $[\lambda_-],[\lambda_+]$ \textbf{bind $S$}, meaning that every complementary region of the union of their supports is a disk or once-punctured disk, then they are the negative and positive directions of a unique Teichm\"uller geodesic which we denote $\tau(\lambda_-,\lambda_+)$.
We note that the binding condition depends only on $|\lambda_-|,|\lambda_+|$, while $\tau(\lambda_-,\lambda_+)$ depends on the projective classes of measures.

We say that a geodesic lamination is \textbf{uniquely ergodic} if it is filling and supports exactly one transverse measure up to scale. By definition, being uniquely ergodic depends only on the support of the lamination.

\subsection{Lengths, intersection numbers, and geodesic currents} \label{lengthsect}

There is a function 
\[
\ell \co \T(S) \times \S \to \mathbb R_+
\]
that assigns a pair $(X, \alpha)$ the length $\ell_X(\alpha)$ of the geodesic representative of $\alpha$ with respect to the hyperbolic metric $X$. This function extends continuously and homogeneously (in the second factor) to a function
\[
\ell \co \T(S) \times \ML(S) \to \mathbb R_+
\] 
called the \textbf{hyperbolic length} function, see \cite{bonahonsurvey}.

There is also a function
\[
\ext \co \T(S) \times \S \to \mathbb R_+
\]
that assigns a pair $(X, \alpha)$ the extremal length $\ext_X(\alpha)$ of the curve $\alpha$ in the Riemann surface $X$---$\ext_X(\alpha)$ is the infimum of the reciprocals of conformal moduli of embedded annuli about $\alpha$.
Kerckhoff proves in \cite{kerckhoff} that this function has a continuous square--homogeneous (in the second factor) extension to $\MF(S)$, namely
\[
\ext \co \T(S) \times \MF(S) \to {\mathbb R}.
\]
Moreover, if $(X,q) \in \Q^*(S)$, then
\[
\ext_X(\F_{\pm}(q)) = \|q\|
\]

This allows the following characterization of the Teichm\"uller metric in terms of extremal length.
\begin{theorem}[Kerckhoff \cite{kerckhoff}]\label{extremalteich} For any $X$ and $Y$ in $\T(S)$
\[
\d_\T(X,Y)=\frac 1 2 \sup \log \left( \frac{\ext_X(\alpha)}{\ext_Y(\alpha)}\right),
\]
where the supremum is taken over all $\alpha$ in $\MF(S)$.
\end{theorem}

Both $\T(S)$ and $\ML(S)$ admit proper embeddings into F. Bonahon's space $\mathfrak{C}(S)$ of \textbf{geodesic currents}: the space of $\pi_1(S)$--invariant Radon measures on the space of geodesics in the universal cover of $S$---we recommend the beautiful \cite{bonahon} for details of what follows.

There is a natural $\mathbb{R}_+$--action on $\mathfrak{C}(S)$ and the quotient $\mathbb{P}\mathfrak{C}(S)$ of $\mathfrak C(S) - \{0\}$ by this action is compact.
The above embeddings descend to embeddings of $\T(S)$ and $\PML(S)$ into $\mathbb{P}\mathfrak{C}(S)$.  There is an ``intersection function'':
\[
\iota \co \mathfrak{C}(S) \times \mathfrak{C}(S) \to {\mathbb R},
\]
which is a continuous symmetric homogeneous (in each factor) function with the following properties.

\begin{theorem}[Bonahon] \label{bony}
Identifying $\T(S)$ and $\ML(S)$ with their images in $\mathfrak{C}(S)$, if $\mu,\nu \in \ML(S)$ and $X \in \T(S)$, then:
\begin{enumerate}
\item $\iota(\mu,\nu) = i(\mu,\nu)$
\item $\iota(X,\mu) = \ell_{X}(\mu)$
\item $\iota(X,X) = \pi^2 |\chi(S)|$
\end{enumerate}
Moreover, $\ML(S)$ consists of precisely those currents $\eta \in \mathfrak{C}(S)$ for which $\iota(\eta,\eta) = 0$.
\end{theorem}

Taking the closure of $\T(S)$ in $\mathbb P \mathfrak C (S)$ provides a compactification of $\T(S)$, as $\mathbb{P}\mathfrak{C}(S)$ is compact.
Properness of the embedding $\T(S) \to \mathfrak{C}(S)$ with part 3 of Theorem \ref{bony} guarantees that any point of $\overline{\T(S)} - \T(S)$ in $\mathbb{P}\mathfrak{C}(S)$ is represented by an element $\eta \in \mathfrak{C}(S)$ satisfying $\iota(\eta,\eta) = 0$.
By the final remark in Theorem \ref{bony}, this is a compactification by $\PML(S)$. Bonahon proves that this is precisely \textbf{Thurston's compactification} \cite{bonahon}, which we write as $\T(S) \cup \PML(S)$.  We comment that this is related to Teichm\"uller's compactification (see the next section) but the two compactifications are different \cite{kerckhoff}.

\subsection{Three theorems of H. Masur}\label{masursect}

We use the following theorems of H. Masur repeatedly.

\begin{theorem}[Criterion for unique ergodicity \cite{masurhaus}]\label{ergcriterion} If the direction of a Teichm\"uller geodesic ray $\tau$ is not uniquely ergodic, then the projection of $\tau$ to the moduli space $\M(S)$ leaves every compact set.
\end{theorem}

\noindent A topological ray $[0,\infty) \to \M(S)$ \textbf{leaves} (or \textbf{exits}) a compact set $\mathcal K$ if the pullback of $\mathcal K$ to $[0,\infty)$ is compact.

\begin{theorem}[Two Boundaries Theorem \cite{twoboundaries}] A Teichm\"uller geodesic ray $\tau$ with direction a uniquely ergodic lamination $\lambda$ converges in $\T(S) \cup \PML(S)$ to $\lambda$. 
\end{theorem}

\noindent So Teichm\"uller's compactification of $\T(S)$ agrees with Thurston's on the set of uniquely ergodic laminations.

\begin{theorem}[Asymptotic Rays Theorem \cite{masuruniquely}]\label{asymptoticrays} Let $X$ and $Y$ be points in $\T(S)$ and let $\sigma$ and $\tau$ be Teichm\"uller geodesic rays from $X$ and $Y$, respectively, with common direction a uniquely ergodic lamination $[\lambda]$ in $\PML(S)$.  Then $\sigma$ and $\tau$ are asymptotic.
\end{theorem}

\noindent Here we say that $\sigma$ and $\tau$ are \textbf{asymptotic} if
\[
\lim_{t \to \infty} \d_\T(\sigma(t),\tau(t)) = 0
\]
for an appropriate choice of unit speed parameterization.

\subsection{Complexes of curves and arcs}\label{complexsect}

Given a surface $Y$ with compact boundary for which the interior $\mathrm{int}(Y)$ is a surface of genus $g$ with $p$ punctures, we let $\xi(Y) = 3g+p$.  We assume throughout that $2 \leq \xi(Y) < \infty$.

A simple closed curve in such a surface $Y$ is \textbf{essential} if it is essential in $\mathrm{int}(Y)$.  A simple (compact) arc is \textbf{essential} if it is homotopically essential relative to $\partial Y$.

Suppose that $\xi(Y) \geq 5$.  Harvey's \textbf{complex of curves} $\C(Y)$ of $Y$ is the simplicial complex whose $k$--cells are collections of isotopy classes of $k+1$ disjoint pairwise non-isotopic essential simple closed curves.

If $\xi(Y) = 4$, then $\mathrm{int}(Y)$ is a sphere with four punctures or a punctured torus.  In these cases, we define $\C(Y)$ to be the graph whose vertices are essential simple closed curves and that two vertices are joined by an edge if they intersect in as few a number of points as is possible for two such curves in $Y$.

When $\xi(Y) \leq 3$, we declare that $\C(Y)$ be empty.

Given a surface $Y$, it is convenient to have a \textbf{complex of arcs} $\mathcal A(Y)$. 
If $\xi(Y) \geq 3$, we define $\mathcal A(Y)$ to be the simplicial complex whose $k$ cells are collections of isotopy classes of $k+1$ disjoint pairwise non-isotopic essential simple closed curves and arcs---where isotopy classes of arcs are defined relative to $\partial Y$.

If $Y$ is a compact annulus, we define $\mathcal A(Y)$ to be the graph whose vertices are isotopy classes of essential arcs  in $Y$ \textit{relative to their endpoints in $\partial Y$} and that two vertices are joined by an edge if they may be realized disjointly.

For any $Y$, we metrize $\C(Y)$ and $\mathcal A(Y)$ by demanding that any simplex is a regular Euclidean simplex with all side lengths equal to one and taking the induced path metric.

We note that when $\xi(Y) \leq 4$ and has no boundary, or when $Y$ is a noncompact annulus, we are declaring $\mathcal A(Y)$ to be empty.

Although the complex $\mathcal A(Y)$ is uncountable when $Y$ is a compact annulus, it is nonetheless quasi-isometric to $\mathbb Z$ \cite{MM2}.

\subsection{Ivanov--Masur--Minsky subsurface projections}\label{subsurfacesect}

A subsurface $Y$ of $S$ is said to be a \textbf{proper domain} if it is not equal to $S$, it is properly embedded (a compact set of $S$ intersects $Y$ in a compact set), and the induced map on fundamental groups is injective.  We further assume that every boundary component of $Y$ is essential (thus the only punctures of $Y$ are also punctures of $S$).
We note that under these assumptions, $\mathcal A(Y)$ is always nonempty.

Fix a hyperbolic metric on $S$.  The definitions for measured laminations which follow are equally valid for their projective classes as well as their supports.  We note that the projections and projection coefficients defined below do not depend on our choice of hyperbolic metric used to realize laminations as geodesic laminations, but the angle $\theta$ and the bound given in Lemma \ref{angle} do.

Given a proper domain $Y$ in $S$, pass to the cover $\widetilde Y$ of $S$ corresponding to the fundamental group of $Y$. Adding the domain of discontinuity for $\pi_1(Y)$ to $\mathbb H^2$ and taking the quotient, we obtain an enlargement $\overline Y$ of $\widetilde Y$ which is homeomorphic to $Y$.

Given $\lambda \in \ML(S)$, we may lift $|\lambda|$ to a (noncompact) geodesic lamination $|\widetilde{\lambda}|$ in $\widetilde Y$.  Compactifying $|\widetilde{\lambda}|$ in $\overline Y$ and identifying any two arcs in the same parallelism class (and disregarding inessential arcs) yields a simplex in $\mathcal A(Y)$---when $Y$ is an annulus, parallelism is defined relative to the endpoints and instead of a simplex, we obtain a set with diameter at most one. This simplex (or set in the annulus case) is the \textbf{projection of $\lambda$ to $Y$}, denoted $\pi_Y(\lambda)$---note that $\pi_Y(\lambda)$ is allowed to be empty.

Given $\mu$ and $\lambda$ in $\ML(S)$, the \textbf{projection coefficient for $\mu$ and $\lambda$ at $Y$} is defined to be
\[
\d_Y(\mu,\lambda) = \mathrm{diam}_{\mathcal A(Y)}(\pi_Y(\mu) \cup \pi_Y(\lambda))
\]
\textit{provided that $\pi_Y(\mu) \neq \emptyset$ and $\pi_Y(\lambda)\neq \emptyset$}.  If either of $\pi_Y(\mu)$ or $\pi_Y(\lambda)$ is empty, we define $\d_Y(\mu,\lambda)=\infty$.

We henceforth write $\mathrm{diam}_Y(\, \cdot \,)$ to denote $\mathrm{diam}_{\mathcal A(Y)}(\, \cdot \,)$.

When $Y$ is an annulus, distance in $\mathcal A(Y)$ is determined by the intersection number: if $\alpha$ and $\beta$ are distinct vertices of $\mathcal A(Y)$, \[ \d_Y(\alpha,\beta) = 1 + i(\alpha,\beta). \]

When convenient, we refer to the core of an annulus $Y$ as $\partial Y$.

Fix a hyperbolic structure on $S$ and $\lambda$ and $\mu$ in $\ML(S)$.
Define the angle
\[ \theta(\mu,\lambda) = \inf_{x \in \mu \cap \lambda} \theta(x,\mu,\lambda) \]
where $\theta(x,\mu,\lambda) \in [0,\frac{\pi}{2}]$ is the smaller of the two angles between tangent lines to $|\mu|$ and $|\lambda|$ at $x$.

\begin{lemma}\label{angle}
Suppose $\alpha \in \S$ is a simple closed geodesic in $S$ and $\mu,\lambda \in \ML(S)$ are two laminations with $\theta = \min\{\theta(\mu,\alpha),\theta(\lambda,\alpha)\} > 0$.
If we let $Y$ denote the annulus with $\partial Y = \alpha$, then we have
\[ \d_Y(\mu,\lambda) \leq 2 \left\lceil 2 \cosh^{-1}(1/\sin(\theta))/\ell(\alpha) \right\rceil + 2. \]
\end{lemma}
\begin{proof}
Let $\widetilde{Y}$ denote the annular cover of $S$ associated to $Y$, $\widetilde{\alpha}$ the lift of $\alpha$ to the core of $\widetilde Y$, and $\widetilde{\mu}$ and $\widetilde{\lambda}$ any lifts of leaves of $|\mu|$ and $|\lambda|$, respectively, that meet $\widetilde{\alpha}$.

Let $\widetilde{\beta}$ denote any geodesic that meets $\widetilde{\alpha}$ orthogonally.
By symmetry and the triangle inequality for $\d_Y$, it suffices to prove that
\[ i(\widetilde{\mu},\widetilde{\beta}) \leq  \left\lceil 2 \cosh^{-1}(1/\sin(\theta)) / \ell(\alpha) \right\rceil. \]

To verify this inequality, further lift the picture to the universal cover ${\mathbb H}^2 \to \widetilde{Y}$.
Let $\widetilde{\alpha}_0$ be a geodesic covering $\widetilde{\alpha}$ that is stabilized by an isometry $\zeta$ generating the cyclic covering group.
Let $\widetilde{\mu}_0$ denote a lift of $\widetilde{\mu}$ intersecting $\widetilde{\alpha}_0$ in a point $x$.
Set $\widetilde{\beta}_0$ to be any lift of $\widetilde{\beta}$ and note that all other lifts of $\widetilde{\beta}$ are of the form $\zeta^n(\widetilde{\beta}_0)$.

Since every point of intersection of $\widetilde \mu$ with $\widetilde \beta$ has exactly one lift on $\widetilde{\mu}_0$, we see that
\[ i(\widetilde{\mu},\widetilde{\beta}) = \sum_{n \in {\mathbb Z}} i(\widetilde{\mu}_0,\zeta^n(\widetilde{\beta}_0)). \]
However, from elementary trigonometric formulae for right triangles we see that a geodesic orthogonal to $\widetilde{\alpha}_0$ will intersect $\widetilde{\mu}_0$ if and only if the distance from this geodesic to $x$ is less than $\cosh^{-1}(1/\sin(\theta(x,\widetilde{\mu}_0,\widetilde{\alpha}_0)))$.
The desired inequality follows from the fact that the translation length of $\zeta$ is $\ell(\alpha)$ and that $\theta \leq \theta(x,\widetilde{\mu}_0,\widetilde{\alpha}_0)$.
\end{proof}

\subsection{Bounded geometry theorems}\label{bddgeomsect}

Minsky's Bounded Geometry Theorem \cite{minskybound} says that a doubly degenerate hyperbolic $3$--manifold homeomorphic to $S \times \mathbb R$ has the length of its shortest geodesic bounded below if and only if the subsurface projection coefficients of its ending laminations are all uniformly bounded above.

K. Rafi has characterized the short curves in hyperbolic structures on a Teichm\"uller geodesic in terms of the intersections of the subsurface projections of its defining laminations \cite{rafi}.  With the global connection between intersection numbers and subsurface projection coefficients described in \cite{MM2}, this yields the following bounded geometry theorem for Teichm\"uller geodesics---the theorem is implicit in the the proof of Theorem 1.5 of \cite{rafi}.

\begin{theorem}[Rafi] \label{rafitheorem}
For every $D > 0$, there exists $\epsilon >0$ such that if $\tau = \tau(\lambda_-, \lambda_+)$ is a Teichm\"uller geodesic with $\lambda_-$ and $\lambda_+$ in $\PML(S)$ satisfying
\[
\d_Y(\lambda_-,\lambda_+) \leq D
\]
for every proper domain $Y \subset S$ with $\xi(Y) \neq 3$, then $\tau$ is $\epsilon$--cobounded.

Conversely, for every $\epsilon  > 0$ there exists $D > 0$ such that if $\tau$ is an $\epsilon$--cobounded Teichm\"uller geodesic with negative and positive directions $\lambda_-$ and $\lambda_+$, then
\[ 
\d_Y(\lambda_-,\lambda_+) \leq D
\]
for every proper domain $Y \subset S$ with $\xi(Y) \neq 3$.
\end{theorem}

\section{Dynamics on $\PML(S)$} \label{boundary}

\subsection{Limit sets} \label{limitsets}
Let $G$ be a subgroup of $\Mod(S)$. The definition of a limit set for the action of $G$ on $\PML(S)$ requires more care than its analogue in the realm of Kleinian groups.
We elaborate here on this notion in our setting.  We primarily follow the notation and conventions of \cite{mccarthypapa}.

A \textbf{weak limit point} for the action of $G$ on $\PML(S)$ is any point $[\lambda] \in \PML(S)$ such that for some $[\mu] \in \PML(S)$ and some infinite sequence of distinct elements $\{ g_n \}_{n=1}^\infty \subset G$, we have $[\lambda] = \lim g_n([\mu])$.  The \textbf{canonical limit set} (for $G$ acting on $\PML(S)$) is the closure of the set of all weak limit points.
A \textbf{limit set} for $G$ is any closed $G$-invariant subset of the canonical limit set.

We say that $G$ is \textbf{dynamically irreducible} if it has a unique non-empty minimal limit set $\PML(S)$.  We call this \textbf{the limit set} and we denote it $\Lambda_G$.  We refer to the points of $\Lambda_G$ as \textbf{limit points} of $G$.
Dynamically irreducible groups fall into two types---see Theorem 4.6 of \cite{mccarthypapa}.  The first type are those which contain a finite index cyclic subgroup (we say it is {\em virtually cyclic}).  Moreover, the finite index cyclic subgroup of $G$ is generated by a pseudo-Anosov mapping class and there is another element of $G$ interchanging the stable and unstable fixed points in $\PML(S)$.  The second (more interesting) type of dynamically irreducible subgroup contains two independent pseudo-Anosov mapping classes.
For a dynamically irreducible group $G$, $\Lambda_G$ can be defined as the closure of the set $\Lambda_0(G)$ of stable laminations of pseudo-Anosov elements of $G$.  Set
\[
Z\Lambda_G = \{[\mu] \in \PML(S)\, |\, i(\mu,\lambda) = 0\ \mbox{for\ some\ } [\lambda] \in \Lambda_G \}.
\]
In \cite{mccarthypapa}, it is shown that $G$ acts properly discontinuously on the set
\[
\Delta_G=\PML(S) - Z\Lambda_G,
\]
called the \textbf{domain of discontinuity} for $G$.

If $G$ is not dynamically irreducible, we say that it is \textbf{dynamically reducible}.
By Theorem 4.6 of \cite{mccarthypapa}, in this case $G$ is either finite, virtually cyclic (virtually) generated by a single pseudo-Anosov mapping class (and contains no element interchanging the stable and unstable fixed points), or is infinite and \textbf{reducible}, which means that there is a nonempty $G$--invariant set $\mathcal{R} \subset \S$ such that for any $\alpha_1,\alpha_2 \in \mathcal{R}$ we have $i(\alpha_1,\alpha_2) = 0$.  We call such a set a \textbf{reduction system} for $G$.
We pause to elaborate on the structure of reducible subgroups---see Chapter 7 of \cite{ivanov}.

If $G$ is infinite and reducible, then there is a \textbf{canonical reduction system} for $G$ characterized as the unique largest reduction system $\mathcal{R}$ with the property that if $\beta \in \S$ is any curve with $i(\alpha,\beta) > 0$ for some $\alpha \in \mathcal{R}$, then there exists $g \in G$ so that $\{ g^k(\beta)\}_{k =1}^\infty$ is an infinite set.  Let $S_1,...,S_n$ denote the components of the complement of the curves of $\mathcal{R}$ in $S$, and we refer to these as the \textbf{components of $G$}.  As $G$ leaves $\mathcal{R}$ invariant, one obtains a homomorphism to the mapping class group of the disjoint union of the components
\[
G \to \Mod(S_1 \sqcup ... \sqcup S_n)
\]
by restriction.
The kernel of the permutation action on the components is a finite index subgroup $G' < G$ leaving each $S_j$ invariant, and we let $G'|_{S_j}$ denote the restriction of $G'$ to $\Mod(S_j)$.
We call a component $S_j$ a \textbf{pseudo-Anosov component} if there is an element $g \in G'$ which is pseudo-Anosov in $G'|_{S_j}$.
A component $S_j$ is called a \textbf{finite component} if the restriction of $G'|_{S_j}$ is finite.
Every component $S_j$ is either pseudo-Anosov or finite---see Theorems 7.11 and 7.14 of \cite{ivanov}.

As dynamically reducible groups do not have unique closed invariant sets on which they act minimally, we make the following declarations of what are to be considered ``\textbf{the limit sets}'' and ``\textbf{the domains of discontinuity}'' of such groups.

If $G$ is finite, we set $\Lambda_G =Z\Lambda_G= \emptyset$ and declare that $\Delta_G = \PML(S)$.

If $G$ is virtually generated by a single pseudo-Anosov mapping class, we define $\Lambda_G$, $Z\Lambda_G$, and $\Delta_G$ as in the dynamically irreducible case.

If $G$ is infinite and reducible, then we follow McCarthy and Papadopoulos and define the limit set and domain of discontinuity as follows.  We let $\mathcal{R}$ denote the canonical reducing system for $G$ and let $S_1,...,S_n$ be the components of $G$, which we number so that for some $m \leq n$, $S_1,...,S_m$ are precisely all the pseudo-Anosov components.  for each $j = 1,...,m$, let $\Lambda_0^j \subset \PML(S)$ denote the set of stable laminations of the pseudo-Anosov elements of $G'|_{S_j}$, considered as laminations in $S$, and $\Lambda^j$ the closure of this set.  The limit set $\Lambda_G$ is defined to be
\[
\Lambda_G = \mathcal{R} \cup \bigcup_{j=1}^m \Lambda^j
\]
The zero set $Z \Lambda_G$ is defined as before to be
\[
Z\Lambda_G = \{[\mu] \in \PML(S)\, |\, i(\mu,\lambda) = 0\ \mbox{for\ some\ } [\lambda] \in \Lambda_G \}.
\]
McCarthy and Papadopoulos prove that the set $\Delta_G = \PML(S) - Z\Lambda_G$ is again a domain of discontinuity for $G$.

\begin{remark} The limit set $\Lambda_G$, its enlargement $Z\Lambda_G$, and the domain $\Delta_G$ have their provenance in Masur's work on the mapping class groups of $3$--dimensional handlebodies \cite{masurhandle}. 
\end{remark}

\subsection{Proper discontinuity revisited} \label{properdiscsect}

The proof that $G$ acts properly discontinuously on $\Delta_G$ given in \cite{mccarthypapa} is easily extended to prove

\begin{theorem} \label{pd}
The action of $G$ on $\T(S) \cup \Delta_G$ is properly discontinuous.
\end{theorem}

The proof follows from a series of lemmata mirroring those in Section 6.2 of \cite{mccarthypapa}.
We invite the reader to visit that paper for further discussion of these ideas.

\begin{remark} McCarthy and Papadopoulos prove that every orbit in $\T(S) \cup \Delta_G$ is discrete---see Section 8 of \cite{mccarthypapa}.
This also follows from Theorem \ref{pd}.
\end{remark}

If $G$ is finite, the theorem is trivial and so we assume that $G$ is infinite for the remainder of this section.

Suppose that $G$ contains a mapping class represented by a pseudo-Anosov homeomorphism $f$ and let $|\L| = \{ \L_-,\L_+ \} \subset \ML(S)$ be unstable and stable measured laminations for $f$ (note that $|\L|$ is a pair of measured lamination, not the support of a lamination).  We comment that these are \textit{measured} laminations, not projective classes---this requires an (arbitrary) choice of representative from the projective classes.  We refer to $|\L|$ as a \textbf{pseudo-Anosov pair} for $f$.  For any $h$ in $G$, we let $|h\L|$ denote the pair $\{ h\L_-,h\L_+ \}$, define
\[
i(\, \cdot \, , |\L|) = \max\{i(\, \cdot \, , \L_-), i(\, \cdot \, ,\L_+)\},
\]
and let
\[
\Delta_{|\L|} = \{ \mu \in \ML(S)\, |\, \forall \, g \in G,\ \  i(\mu,|\L|) \leq i(\mu,|g\L|)\}.
\]
This set is $\mathbb R_+$--invariant and so defines a subset of $\PML(S)$ which we also call $\Delta_{|\L|}$.

McCarthy and Papadopoulos show that $\Delta_{|\L|}'=\Delta_{|\L|} \cap \Delta_G$ is a fundamental domain for the action of $G$ on $\Delta_G$.

\begin{remark} If a group $\Gamma$ acts on a topological space $\mathcal X$, we say that a closed subset $\mathcal D \subset \mathcal X$ is a \textbf{fundamental domain} for the action if $\{ \gamma \mathcal D\, | \, \gamma \in \Gamma \}$ is a locally finite covering of $\mathcal X$.
\end{remark}

We extend the function $i(\, \cdot \, , |\L|)$ defined on $\ML(S)$ in the obvious way to a function $\iota(\, \cdot \, , |\L|)$ defined on the union $\ML(S) \cup \T(S)$ in $\mathfrak{C}(S)$.
We define
\[
\widehat{\Delta}_{|\L|} = \{ X \in \T(S) \cup \ML(S) \, |\, \forall \, g \in G,\ \  \iota(X,|\L|) \leq \iota(X,|g\L|)\}.
\]
Again, this set is $\mathbb R_+$--invariant and so defines a set in $\T(S) \cup \PML(S)$ that we also call $\widehat{\Delta}_{|\L|}$.
As with $\Delta_{|\L|}$, one readily checks that $\widehat{\Delta}_{|\L|}$ is closed.

As in \cite{mccarthypapa}, we note that for any $g \in G$, we have
\[
g(\widehat{\Delta}_{|\L|}) = \widehat{\Delta}_{|g\L|}.
\]

As in the proof of proper discontinuity on $\Delta_G$ given in \cite{mccarthypapa}, special attention must be paid when $G$ is reducible. In this case, we proceed as follows---see Section \ref{limitsets} for notation.
For each component $S_j$ of $G$, let $|\L^j| = \{ \L^j_-,\L^j_+ \}$ be a pseudo-Anosov pair for some pseudo-Anosov automorphism in $G'|_{S_j}$, viewed as laminations in $\ML(S)$.
We let $\upsilon$ denote any curve that non-trivially intersects each component of $\mathcal{R}$, and let $|\L|$ denote the union of $\upsilon$ and each $|\L^j|$ for $j = 1,...,m$.
$|\L|$ is called a \textbf{complete system for $G$}.
We define 
\[
\iota(\cdot,|\L|) = \max\{\iota(\cdot,\L^i_-), \iota(\cdot, \L^i_+), \iota(\cdot, \upsilon) \, | \, 1\leq i \leq n \}
\]
and define the sets $\Delta_{|\L|}$, $\Delta_{|\L|}'$, and $\widehat{\Delta}_{|\mathcal{L}|}$ exactly as before.

The first lemma we need is the following (compare Lemma 6.11 \cite{mccarthypapa}).  In the following discussion, points of $\T(S) \cup \PML(S)$ will be enclosed in brackets $[X]$, and we will remove the brackets $X$ to denote a representative of $[X]$ in $\T(S) \cup \ML(S)$.  This is only relevant when $[X] \in \PML(S)$, in which case $X$ is a representative of the projective class $[X]$.  For $[X] \in \T(S)$, we have $X = [X]$.

\begin{lemma} \label{611}
Let $[X] \in \T(S) \cup \Delta_G$, and $\{g_n\}$ be an infinite sequence of distinct mapping classes in $G$.
Then the sequence of numbers $\{\iota(X,|g_n\L|)\}$ is unbounded.
\end{lemma}
\begin{proof} Suppose that $G$ is irreducible.

Upon passing to a subsequence, the hypothesis implies that one of the sequences $\{ g_n\L_- \}$ or $\{g_n\L_+ \}$ is unbounded in $\ML(S)$ (see Lemma 2.6 of \cite{mccarthypapa}).  That is, there is a curve $\alpha$ so that one of the sequences $\{ i(\alpha,g_n \L_+) \}$ or $\{i(\alpha,g_n \L_-)\}$ is tending to infinity.
Without loss of generality, we assume $\{ g_n\L_+ \}$ is unbounded.
So there is a sequence of positive real numbers $\{ r_n \}$ tending to $0$ such that $r_ng_n\L_+ \to \mu$ in $\ML(S)$.
Since $[\L_+] \in \Lambda_G$, so is $[\mu]$.
We have
\[
\iota(X,r_ng_n\L_+) = r_n\iota(X,g_n \L_+) \leq r_n \iota(X,|g_n\L|).
\]
If $\iota(X,|g_n\L|)$ were bounded independent of $n$, the numbers $\iota(X,r_ng_n\L_+)$ would converge to zero, implying that $\iota(X,\mu) = 0$.
For $[X] \in \Delta_G = \PML(S) - Z\Lambda_G$ this is an obvious contradiction.
If $[X] \in \T(S)$, this would mean that $\ell_X(\mu) = 0$, which is also impossible.

If $G$ is reducible, one of $\{g_n \L^i_-\}$, $\{g_n \L^i_+\}$, and $\{g_n \upsilon\}$ is unbounded in $\ML(S)$ (see Lemma 7.6 of \cite{mccarthypapa}), and the proof continues as in the irreducible case.
\end{proof}

The next fact we need is our version of Proposition 6.13 of \cite{mccarthypapa}.

\begin{lemma} \label{613}
For every $[X] \in \T(S) \cup \Delta_G$, there exists $g \in G$ such that $[X] \in g \widehat{\Delta}_{|\L|}$. 
\end{lemma}
\begin{proof}
By Lemma \ref{611} the set $\{\iota(X,|g\L|)\}_{g \in G}$ has no infinite bounded subsets.
It follows that there is a minimum $\iota(X,|g\L|)$ for some $g \in G$ and hence $[X] \in \widehat{\Delta}_{|g\L|} = g \widehat{\Delta}_{|\L|}$.\end{proof}

We now turn to the analog of Proposition 6.14 of \cite{mccarthypapa}.

\begin{lemma} \label{614}
Let $\mathcal K$ be a compact set in $\T(S) \cup \Delta_G$.
Then the set of mapping classes $\{g \in G \, | \, \mathcal K \cap g \widehat{\Delta}_{|\L|} \neq \emptyset \}$ is finite.
\end{lemma}
\begin{proof}
Suppose there is an infinite sequence $\{g_n\}$ of distinct elements of $G$ such that $\mathcal K \cap g_n \widehat{\Delta}_{|\L|} = \mathcal K \cap \widehat{\Delta}_{|g_n\L|} \neq \emptyset$ for every $n$ and let $[X_n] \in \mathcal K \cap \widehat{\Delta}_{|g_n\L|}$.
In particular, we have $\iota(X_n,|g_n\L|) \leq \iota(X_n,|\L|)$.

Suppose that $G$ is irreducible.

As above (and in \cite{mccarthypapa}), one of the sequences $\{g_n\L_+\}$ or $\{g_n\L_-\}$ is unbounded, and we assume without loss of generality that it is the first.
After passing to subsequences, there is a pair of sequences $\{r_n\}$ and $\{t_n\}$ of positive real numbers, the first tending to $0$, so that
\[
r_ng_n\L_+ \to \mu \in \ML(S) \quad \mbox{ and } \quad t_nX_n \to X \in \mathfrak{C}(S).
\]

As in \cite{mccarthypapa}, it follows from continuity of $\iota$ that
\begin{align*}
\iota(X,\mu) & = \lim \iota(t_nX_n,r_ng_n\L_+)\\
         & \leq \lim r_n \iota(t_nX_n,|g_n\L|)\\
         & \leq \lim r_n \iota(t_nX_n,|\L|)\\
         & = (\lim r_n)(\lim \iota(t_nX_n,|\L|))\\
         & = (0)(\iota(X,|\L|))\\
         & = 0.
\end{align*}
As in the proof of Lemma \ref{611}, this contradicts the fact that $\mu \in \Lambda_G$ and $[X] \in \mathcal K \subset \T(S) \cup \Delta_G$.

In the reducible case, one of $\{g_n \L^i_-\}$, $\{g_n \L^i_+\}$, and $\{g_n \upsilon\}$ is unbounded and again the proof is formally identical to the irreducible case.
\end{proof}

We may now prove Theorem \ref{pd} (compare to the proof of Theorem 6.16 of \cite{mccarthypapa}).
\begin{proof}
Let $\mathcal K \subset \T(S) \cup \Delta_G$ be compact. We show that the set $\{g \in G\, |\, g\mathcal K \cap \mathcal K \neq \emptyset \}$ is finite.

By Lemma \ref{614} the set
\[
\{ g \in G \, | \, \mathcal K \cap g \widehat{\Delta}_{|\L|} \neq \emptyset \}
\]
is finite, and we name its elements $g_1, \, \ldots \,  ,g_N$.
With Lemma \ref{613}, we see that
\[
\mathcal K \subset \bigcup_{j=1}^N g_j \widehat{\Delta}_{|\L|} \, .
\]

Now, if $g \mathcal K \cap \mathcal K \neq \emptyset$, then $g(g_j \widehat{\Delta}_{|\L|}) \cap \mathcal K \neq \emptyset$ for some $j \in \{ 1, \, \ldots \,  ,N \}$.
Since 
\[
g(g_j \widehat{\Delta}_{|\L|}) = (gg_j)\widehat{\Delta}_{|\L|},
\]
it follows that $gg_j = g_i$ for some $i \in \{ 1, \, \ldots \, ,N \}$.
In particular,
\[
\{g \in G \, | \, g\mathcal K \cap \mathcal K \neq \emptyset \} \subset \{g_ig_j^{-1}\}_{i,j=1}^N.
\]
Since the set on the right is finite, so is the one on the left.
\end{proof}

We have also established

\begin{proposition} \label{fundomain}
Let $G$ be an infinite subgroup of $\Mod(S)$. Then 
\[
\widehat{\Delta}'_{|\L|} = \widehat{\Delta}_{|\L|} \cap (\T(S) \cup \Delta_G)
\]
is a fundamental domain for the action of $G$ on $\T(S) \cup \Delta_G$. \qed
\end{proposition}

\subsection{Conical limit points and compact fundamental domains} \label{sectionconicallimitpoints}

A point $[\lambda]$ in $\Lambda_G$ is a \textbf{conical limit point} if for any Teichm\"uller geodesic ray $\tau$ with direction $[\lambda]$, there is a number $R > 0$ such that some $G$--orbit intersects the $R$--neighborhood of $\tau$ in an infinite set. Note that a conical limit point is uniquely ergodic by Masur's criterion---as the projection of the geodesic $\tau$ to the moduli space $\M(S)$ must return to a bounded neighborhood of a point infinitely often and so cannot leave every compact set. In particular, any geodesic ray $\tau$ whose direction is a conical limit point in fact terminates at that point in Thurston's compactification of $\T(S)$, by Masur's Two Boundaries Theorem.

If $\tau$ is a Teichm\"uller geodesic ray emanating from a point $X$ with direction $[\lambda]$ whose $R$--neighborhood contains infinitely many points of a $G$--orbit, then any geodesic ray $\sigma$ terminating at $[\lambda]$ has an $R'$--neighborhood containing infinitely many points from that orbit.  To see this, note that $\sigma$ and $\tau$ are asymptotic by Masur's Asymptotic Rays Theorem.  In particular, $\sigma$ and $\tau$ are at a finite Hausdorff distance $A$ from each other and it suffices to take $R' = R + A$. So, to verify that a limit point is conical, we need only consider a single ray.


\begin{theorem}\label{ccconical} If $G$ is a convex cocompact subgroup of $\Mod(S)$, then every one of its limit points is conical. 
\end{theorem}
\begin{proof} Let $[\lambda] \in \Lambda_G$ and let $\tau$ be a geodesic ray in $\WH_G$ with direction $[\lambda]$.  Since $G$ acts cocompactly on $\WH_G$, there is a positive number $R$ such that every point of the image of $\tau$ in $\WH_G/G$ is a distance at most $R$ from any fixed point $X_0$ in $\WH_G/G$.  So, if $X$ is a point in the preimage of $X_0$, $\tau$ stays within $R$ of $G X$.
\end{proof}


For the remainder of this subsection, it is convenient to switch points of view and work primarily with measured foliations instead of measured laminations.  As such, we let $|\F| = \{ \F_- , \F_+ \}$ denote the foliations associated to a pseudo-Anosov pair $\{ \L_- , \L_+ \}$.  If $\lambda \in \ML(S)$ is associated to $F \in \MF(S)$, then we have
$i(F,|\F|) = i(\lambda,|\L|)$.  Also, if $g \in G$, then we write $|g \F| = \{g \F_-,g \F_+ \}$.

The pair $|\F|$ determines a unique point $(X,q) = (X_{|\F|},q_{|\F|}) \in \Q^*(S)$ with the property that $\F_- = c \F_-(q_{|\F|})$ and $\F_+ = c \F_+(q_{|\F|})$ for some $c > 0$.  Scaling $|\F|$, we may assume that $(X,q) \in \Q^1(S)$ (and so $c = 1$).

As in Section 3 of \cite{minskycrelle}, we see that for any $Y \in \T(S)$ and any $\alpha,\beta \in \MF(S)$
\[
\ext_Y(\alpha)\, \ext_Y(\beta) \geq i(\alpha,\beta)^2.
\]
Since $\F_-$ and $\F_+$ are the horizontal and vertical foliation of $q$, and since $\|q\| = 1$, we see that
\[
\ext_X(\F_-) = \ext_X(\F_+) = 1.
\]
Thus, for any $F \in \MF(S)$ with $\lambda \in \ML(S)$ the associated lamination, we have
\begin{equation} \label{extremal}
\ext_X(F) = \ext_X(F)\, \ext_X(\F_{\pm}) \geq i(F,|\F|)^2 = i(\lambda,|\L|)^2.
\end{equation}


\begin{theorem}\label{conicaldelta} Let $[\lambda]$ be a conical limit point of $G$.  Then $[\lambda] \notin \widehat{\Delta}_{|\L|}$.
\end{theorem}
\begin{proof} Let $|\F|$ and $(X,q) = (X_{|\F|},q_{|\F|}) \in \Q^1(S)$ be as above. Let $[\lambda]$ be a conical limit point of $G$ and $F$ the measured foliation associated to $\lambda$ (an arbitrary choice from the projective class of $[\lambda]$). Let $\tau$ be the Teichm\"uller geodesic emanating from $X$ and terminating at $[F]$ defined by a unit norm quadratic differential $\omega$ at $X$.  Thus, choosing the representative of the projective class appropriately, we may assume that $F$ is the vertical foliation of $\omega$.  Note also that for $g$ in $G$, $g(X_{|\F|}) = X_{|g \F|}$.

Now, by the conical hypothesis there is a positive number $K$ and an infinite set $\{g_n\}_{n=0}^\infty \subset G$ such that $g_n(X)$ is a distance at most $\frac 1 2 \log K$ from a point $Y_n$ on $\tau$. 
Since $F$ is the vertical foliation of $\omega$, we have
\[
\lim_{n \to \infty} \ext_{Y_n}(F) = 0.
\]
Since $\d_\T(Y_n, g_n(X))$ is no more than $\frac 1 2 \log K$,
\[
\ext_{g_n(X)}(\alpha) \leq K \, \ext_{Y_n}(\alpha)
\]
for all $\alpha$ in $\MF(S)$ by Theorem \ref{extremalteich} and so
\[
\lim_{n \to \infty} \ext_{g_n(X)}(F) = 0.
\]
We conclude that
\[
\lim_{n \to \infty} i(\lambda, |g_n \L|) = \lim_{n \to \infty} i(F, |g_n \F|) = 0
\]
by (\ref{extremal}) and so $[\lambda] \notin \widehat{\Delta}_{|\L|}$.
\end{proof}


\begin{theorem}\label{conicalcompact} Let $G$ be a subgroup of $\Mod(S)$ such that every point in $\Lambda_G$ is conical.  Then there is a compact fundamental domain for the action of $G$ on $\T(S) \cup \Delta_G$.
In particular, convex cocompact groups act cocompactly on $\T(S) \cup \Delta_G$.
\end{theorem}
\begin{proof} By Proposition \ref{fundomain}, the set $\widehat{\Delta}_{|\L|}'$ is a fundamental domain for the action of $G$ on $\T(S) \cup \Delta_G$. Since conical limit points are uniquely ergodic, every lamination in $\Lambda_G$ is uniquely ergodic and so $Z\Lambda_G = \Lambda_G$.  In particular, $\Delta_G = \PML(S) - \Lambda_G$. By Theorem \ref{conicaldelta}, $\widehat{\Delta}_{|\L|}' = \widehat{\Delta}_{|\L|}$. But the set $\widehat{\Delta}_{|\L|}$ is a closed subset of $\T(S) \cup \PML(S)$, and is thus compact.
\end{proof}

\subsection{The weak hull} \label{weakhullssection}


Let $\A$ be a closed subset of $\PML(S)$.  If $\A$ has the property that for every $[\lambda_-] \in \A$, there exists a $[\lambda_+] \in \A$ such that $[\lambda_-]$ and $[\lambda_+]$ bind $S$, then we define the \textbf{weak hull} $\WH_\A$ of $\A$ to be the union of all geodesics $\tau(\lambda_-, \lambda_+)$ in $\T(S)$ with $[\lambda_-]$ and $[\lambda_+]$ elements of $\A$ that bind.  If $\A$ does not have this property then we say that the weak hull is not defined.  A set $\WH$ is a \textbf{weak hull} if it is $\WH_\A$ for some closed $\A \subset \PML(S)$ with the aforementioned property.

Note that if $G$ is an infinite irreducible subgroup of $\Mod(S)$, then $\Lambda_G$ possesses a nonempty weak hull $\WH_G = \WH_{\Lambda_G}$:  by a theorem of Ivanov (\cite{ivanov}, Corollary 7.14), there exists a pseudo-Anosov automorphism in $G$ with stable lamination $[\lambda]$, and any other lamination in $\Lambda_G$ will bind with $[\lambda]$.


\subsection{Compact fundamental domains in $\Delta_G$ cobound the hull}

Having a compact fundamental domain for the action on the domain of discontinuity is a restrictive condition in itself, and in particular suffices to cobound the hull.

\begin{theorem}\label{domainsub} Let $G$ be a subgroup of $\Mod(S)$. If $\Delta_G \neq \emptyset$ and $G$ acts cocompactly on $\Delta_G$, then every lamination in $\Lambda_G$ is uniquely ergodic, $Z\Lambda_G = \Lambda_G$, and $\WH_G$ is defined and cobounded. Furthermore, $G$ has a finite index subgroup all of whose non-identity elements are pseudo-Anosov.
\end{theorem}

The first step to prove this theorem is to prove that every lamination in $\Lambda_G$ is filling.
To this end, we define
\[
Z Z \Lambda_G = \{[\beta] \in \PML(S)\, |\, i(\beta,\mu) = 0\ \mbox{for\ some\ } [\mu] \in Z\Lambda_G \}.
\]
and prove

\begin{theorem}[Insomnia]\label{zz} Let $G$ be a subgroup of $\Mod(S)$. If $\Delta_G \neq \emptyset$, and $G$ acts cocompactly on $\Delta_G$, then 
\[
Z Z \Lambda_G = Z \Lambda_G.
\]
Moreover, every lamination in $\Lambda_G$ is filling.
\end{theorem}
\begin{proof} If $G$ is finite, both $Z\Lambda_G$ and $ZZ\Lambda_G$ are empty, and so we assume that $G$ is infinite.
We begin by proving the last statement, assuming the first.

If there is a non-filling lamination in $\Lambda_G$, then there is a (projective class of laminations supported on a) simple closed curve in $Z\Lambda_G$.
Let $V(Z\Lambda_G) \subset \C(S)$ denote the set of all simple closed curves in $\C(S)$ that lie in $Z\Lambda_G$.
Since $\Delta_G$ is an open set, it contains a simple closed curve that is not in $Z\Lambda_G$, and hence $\C(S) \neq V(Z\Lambda_G) \neq \emptyset$.
Since $\C(S)$ is connected, there is a simple closed curve $[\alpha]$ at a distance $1$ from $V(Z\Lambda_G)$.
The curve $[\alpha]$ is thus disjoint from some element of $Z\Lambda_G$, and is not in $Z\Lambda_G$ (as it is a positive distance from $V(Z\Lambda_G)$).
That is, $[\alpha] \in ZZ\Lambda_G - Z\Lambda_G$, contradicting the first part of the theorem.

We now proceed to the proof of the first statement.
If $G$ is irreducible, let $|\L| = \{\L_-,\L_+\}$ be a pseudo-Anosov pair for $G$. If $G$ is reducible, let $|\L|$ be a complete system for $G$ (see Section \ref{properdiscsect}). Let $\mathcal K$ be a compact fundamental domain for the action of $G$ on $\Delta_G$. By Propositions 6.14 and 7.10 of \cite{mccarthypapa}, the set of mapping classes
\[
\{g \in G\, |\, \mathcal K \cap \Delta_{|g\L|}' \neq \emptyset \}
\]
is finite.  Since
\begin{align*}
|\{g \in G\, |\, \mathcal K \cap \Delta_{|g\L|}' \neq \emptyset \}| 
& = |\{g \in G\, |\, g^{-1} \mathcal K \cap g^{-1} \Delta_{|g\L|}' \neq \emptyset \}|\\
& = |\{g \in G\, |\, g^{-1} \mathcal K \cap \Delta_{|\L|}' \neq \emptyset \}|\\
& = |\{g \in G\, |\, g \mathcal K \cap \Delta_{|\L|}' \neq \emptyset \}|
\end{align*}
and $\mathcal K$ is a fundamental domain, we conclude that $\Delta_{|\L|}'$ is compact.

Suppose to the contrary that there is a lamination $[\beta]$ in $Z Z \Lambda_G - Z\Lambda_G$. So $i(\beta,\lambda) \neq 0$ for all $[\lambda]$ in $\Lambda_G$, and there is a $[\mu]$ in $Z\Lambda_G$ such that $i(\beta,\mu) = 0$. Note that $[\mu]$ cannot be filling, lest $[\beta]$ be an element of $Z\Lambda_G$. So, in fact, there is a simple closed curve $[\alpha]$ in $Z\Lambda_G$ such that $i(\beta, \alpha)=0$. To see this, take $[\alpha]$ be a component of the boundary of the smallest $\pi_1$--injective subsurface containing $|\mu|$.  If $\lambda$ is any lamination with $i(\lambda,\mu) = 0$, then notice that one also has $i(\lambda,\mu) = 0$.

For $t \in [0,1]$, let $\nu_t = (1-t) \alpha + t \beta$.  That is, for $t \in (0,1)$, $\nu_t$ is supported on $|\alpha| \cup |\beta|$ and assigns to each arc transverse to $|\alpha| \cup |\beta|$ the sum of the transverse measures for $\alpha$ and $\beta$ weighted by $(1-t)$ and $t$, respectively.  We also have $\nu_0 = \alpha$ and $\nu_1 = \beta$.
For any interval ${\mathbb J} \subset [0,1]$, write
\[
\nu_{\mathbb J} = \{ \nu_t\, |\, t \in \mathbb J \}
\]
and let $[\nu_{\mathbb J}]$ denote the image in $\PML(S)$.
Since $i(\beta,\lambda) \neq 0$ for every $[\lambda] \in \Lambda_G$, the entire interval $[\nu_{(0,1]}]$ is contained in $\Delta_G$.

Now, for any number $C$, the set
\[
\{ g \in G\, |\, i(\beta,|g\L|) \leq C \}
\] 
is finite, by Lemmata 6.11 and 7.7 of \cite{mccarthypapa}. 

When $G$ is irreducible, the laminations $\L_-$ and $\L_+$ are associated to measured foliations $\F_-$ and $\F_+$. As in Section \ref{sectionconicallimitpoints}, these determine a point $(X,q) = (X_{|\F|},q_{|\F|}) \in \Q^1(S)$ with $\F_-(q) = \F_-$ and $\F_+(q) = \F_+$.

For any simple closed curve $\gamma$, we have
\[
\frac{1}{\sqrt{2}} \ell_q(\gamma) \leq i(\gamma, |\F|) = i(\gamma,|\L|) \leq \ell_q(\gamma)
\]
where $\ell_q\co \S \to \mathbb R_+$ is the function that assigns a curve its $q$--length. Now, for any constant $C$, the set
\[
\{g\alpha\, |\, i(g\alpha, |\L|) \leq C \}
\]
is finite, since the length spectrum of the $q$--metric is discrete. 

When $G$ is reducible, the set 
\[
\{g\alpha\, |\, i(g\alpha, |\L|) \leq C \}
\]
is again finite.
To see this, consider the analogous quadratic differentials $q_j$ on $S_j$ determined by $\L^j_-$ and $\L^j_+$.
The same comparison of length and intersection number shows that the set is finite up to Dehn twisting along the components of $\mathcal{R}$.
However, an infinite collection of curves that differ only by twists in the components of $\mathcal{R}$ will have unbounded intersection numbers with $\upsilon$, and therefore the intersection numbers with $|\L|$ will be unbounded.

Note that since $\Delta_{|\L|}'$ is a compact fundamental domain and $[\nu_{(0,1]}]$ is closed in $\Delta_G$ and non-compact, the set
\[
\{ g \in G \, |\, [\nu_{(0,1]}] \cap \Delta_{|g\L|}' \neq \emptyset \}
\]
is infinite.

For any $T \in (0,1]$, the set $[\nu_{[T,1]}]$ is compact and contained in $\Delta_G$. By Propositions 6.14 and 7.10 of \cite{mccarthypapa}, $[\nu_{[T,1]}]$ only intersects finitely many translates of $\Delta_{|\L|}'$.  For any $t\in (0,1]$, $[\nu_t] \in h\Delta_{|\L|}'$ for some $h$ in $G$, since $\Delta_{|\L|}'$ is a fundamental domain for the action of $G$ on $\Delta_G$. So, we may choose a sequence $t_n$ tending to zero such that $[\nu_{t_n}] \in g_n \Delta_{|\L|}'$ and $\{g_n\}$ is an infinite set---and we do so. By the definition of $\Delta_{|g_n\L|}$, we have
\begin{align*}
i(\nu_{t_n},|g_n\L|) & \leq i(\nu_{t_n}, |\L|)\\
& = (1-t_n)\, i(\alpha, |\L|) + t_n\, i(\beta,|\L|)
\end{align*}
and so the $i(\nu_{t_n},|g_n\L|)$ are uniformly bounded by some number $r$.  In particular,
\begin{align*}
i(\nu_{t_n},|g_n\L|) &  = (1-t_n)\, i(\alpha, |g_n\L|) + t_n\, i(\beta,|g_n\L|)\\
& \leq r
\end{align*}
and so
\begin{align*}
i(g_n^{-1} \alpha, |\L|) & = i(\alpha,|g_n\L|)\\
& \leq r/(1-t_n)\\
& \leq 2r
\end{align*}
when $n$ is large. We conclude that the numbers $i(g_n^{-1} \alpha, |\L|)$ are all bounded by some number $R$. Since the set
\[
\{g\alpha\, |\, i(g\alpha, |\L|) \leq R \}
\]
is finite, we may pass to a subsequence so that 
\[
\{g_n^{-1}\alpha\}_{n\in \mathbb N} = \{\alpha'\}
\]
for some simple closed curve $\alpha'$, and we do so. 

Again, by the definition of $\Delta_{|g_n\L|}$, we have
\begin{equation}\label{zzeq}
i(\nu_{t_n},|g_n\L|) \leq i(\nu_{t_n},|g\L|)
\end{equation}
for all $g$ in $G$. 
Since $\{g_n\}$ is infinite and $[\beta] \in \Delta_G$, the sequence of numbers
\[
i(g_n^{-1}\beta,|\L|) = i(\beta,|g_n\L|)
\]
is unbounded, by Lemmata 6.11 and 7.7 of \cite{mccarthypapa}. So, for some $n$, we have
\begin{align*}
i(\nu_{t_n}, |g_n \L|)
&  = (1-t_n)\, i(g_n^{-1}\alpha, |\L|) + t_n\, i(g_n^{-1}\beta,|\L|)\\
&  = (1-t_n)\, i(\alpha', |\L|) + t_n\, i(g_n^{-1}\beta,|\L|)\\
&  > (1-t_n)\, i(\alpha', |\L|) + t_n\, i(g_1^{-1}\beta,|\L|)\\
&  = (1-t_n)\, i(g_1^{-1}\alpha, |\L|) + t_n\, i(g_1^{-1}\beta,|\L|)\\
&  = i(\nu_{t_n},|g_1\L|)
\end{align*}
which contradicts (\ref{zzeq}).
\end{proof}


With Theorem \ref{zz} in hand, we prove Theorem \ref{domainsub}.

\begin{proof}[Proof of Theorem \ref{domainsub}]
We begin by proving that there exists a constant $D > 0$ such that for any two distinct points $[\lambda_-],[\lambda_+] \in \Lambda_G$ and any proper domain $Y \subset S$ with $\xi(Y) \neq 3$, we have
\[
\d_Y(\lambda_-, \lambda_+) \leq D.
\]
We then appeal to Rafi's Theorem and Masur's criterion for unique ergodicity to see that every lamination in $\Lambda_G$ is uniquely ergodic, $Z\Lambda_G = \Lambda_G$, and the weak hull $\WH_G$ of $\Lambda_G$ is defined and cobounded.
Our proof is a modification of the proof given in \cite{minskyinj} in the case when $G$ is cyclic.

Let $\mathcal K \subset \Delta_G$ be a compact set.  Fix a hyperbolic metric on $S$ and use this to realize $\ML(S)$.

For $[\lambda]$ in $\Lambda_G$ and $[\kappa]$ in $\mathcal K$, let $L([\lambda], [\kappa])$ denote the supremum of lengths of arcs of $|\lambda| \cap (S-|\kappa|)$.  This is finite since $\lambda$ is filling (by Theorem \ref{zz}). We claim that $L(\, \cdot \, , \, \cdot \,)$ is bounded on $\Lambda_G \times \mathcal K$.   

Suppose to the contrary that there are sequences $[\lambda_i]$ in $\Lambda_G$ and $[\kappa_i]$ in $\mathcal K$ such that $L([\lambda_i], [\kappa_i])$ tends to infinity with $i$.  Since $\Lambda_G$ and $\mathcal K$ are compact, we may assume that the $[\lambda_i]$ tend to a lamination $[\lambda]$ in $\Lambda_G$, the $[\kappa_i]$ to a lamination $[\kappa]$ in $\mathcal K$. Since the $L([\lambda_i], [\kappa_i])$ are tending to infinity, we have a sequence of geodesic arcs $\alpha_i$  in $|\lambda_i| \cap (S-|\kappa_i|)$ whose lengths are tending to infinity. The Hausdorff limit of the $\alpha_i$ is a diagonal extension of $|\lambda|$; we conclude that $|\kappa|$ does not transversely intersect that extension, and hence $\kappa$ has zero intersection number with $\lambda$.  But this means that $[\kappa]$ is an element of $Z \Lambda_G$, contradicting the fact that $\mathcal K \cap Z\Lambda_G$ is empty.

So, there is a constant $L = L(G,\mathcal K)$ such that for any $[\lambda]$ in $\Lambda_G$ and any $[\kappa]$ in $\mathcal K$, the length of any arc in $|\lambda| \cap (S - |\kappa|)$ is bounded above by $L$.

A similar argument shows that there is a constant $\Theta = \Theta(G,\mathcal K) > 0$ such that for all $[\lambda]$ in $\Lambda_G$ and $[\kappa]$ in $\mathcal K$, the angle $\theta([\lambda],[\kappa])$ between leaves of $|\lambda|$ and those of $|\kappa|$ is at least $\Theta$.

Note that if $Y$ is a proper domain that is not an annulus and $\lambda$ is a geodesic lamination, the projection $\pi_Y(\lambda)$ may be obtained by realizing the boundary components of $Y$ as geodesics and intersecting $|\lambda|$ with $\mathrm{int}(Y)$. Therefore, by the Keen--Halpern Collar Lemma \cite{keen,halpern}, there is a constant $M=M(G,\mathcal K)$ such that for any pair $[\lambda_-]$ and $[\lambda_+]$ in $\Lambda_G$ and any proper subdomain $Y$ with $\xi(Y) \geq 4$ and $\partial Y$ in $\mathcal K$,
\[
i(\pi_Y(\lambda_-), \pi_Y(\lambda_+)) \leq M.
\]
This implies the existence of a constant $B=B(G,\mathcal K)$ such that
\[
\d_Y(\lambda_-, \lambda_+) \leq B
\]
whenever $Y$ is not an annulus and $\partial Y$ is an element of $\mathcal K$.

When $Y$ is an annulus, 
\[ 
\d_Y(\lambda_-,\lambda_+) \leq 2 \left\lceil 2 \cosh^{-1}(1/\sin(\Theta))/\ell(\partial Y) \right\rceil + 2 
\]
whenever $\partial Y$ is in $\mathcal K$, by Lemma \ref{angle}. The injectivity radius of our chosen hyperbolic metric bounds $\ell(\partial Y)$ from below, and so there is a constant $C=C(G,\mathcal K)$ such that
\[
\d_Y(\lambda_-, \lambda_+) \leq C
\]
whenever $Y$ is an annulus and $\partial Y$ lies in $\mathcal K$.

Letting $D = D(G,\mathcal K) = \max\{B,C \}$, we conclude that 
\[
\d_Y(\lambda_-, \lambda_+) \leq D
\]
whenever $Y$ is a proper domain with $\partial Y$ in $\mathcal K$.\\

If there is a compact fundamental domain $\mathcal K$ for the action of $G$ on $\Delta_G$, we have the bound
\[
\d_Y(\lambda_-, \lambda_+) \leq D
\]
for all $Y$ with $\partial Y$ in $\Delta_G$ and all pairs $[\lambda_-]$ and $[\lambda_+]$ in $\Lambda_G$, since $\Lambda_G$ is $G$--invariant and $[\lambda_-]$ and $[\lambda_+]$ were arbitrary.  

As the laminations in $\Lambda_G$ are filling by Theorem \ref{zz}, given a proper domain $Y$, $\partial Y$ is an element of $\Delta_G$, and we have the desired bound for all proper subdomains.

By Rafi's Theorem, the geodesics joining distinct points in $\Lambda_G$ are uniformly co-bounded.

Let $[\lambda]$ be an element of $\Lambda_G$.  We may find a $[\lambda']$ in $\Lambda_G$ such that $\lambda$ and $\lambda'$ bind $S$.  To see this, first note that $G$ is either finite, in which case the conclusions of the theorem are trivial, or $G$ contains a pseudo-Anosov mapping class; for if not, then $Z\Lambda_G$ would contain a simple closed curve, which is prohibited by the filling hypothesis. Now, let $[\lambda_-]$ and $[\lambda_+]$ be the unstable and stable laminations of a pseudo-Anosov mapping class $g$ in $G$. If $[\lambda]$ is an element of $\{[\lambda_-],[\lambda_+]\}$, we are done.  If not, $[\lambda]$ and $[\lambda_-]$ bind $S$. In any case, there is a Teichm\"uller geodesic with directions $[\lambda]$ and $[\lambda']$. This geodesic is cobounded by the previous paragraph and Masur's criterion tells us that $[\lambda]$ is uniquely ergodic. We conclude that $Z\Lambda_G = \Lambda_G$.

Since every lamination in $\Lambda_G$ is uniquely ergodic, every pair of points in $\Lambda_G$ are joined by a Teichm\"uller geodesic.  So the weak hull of $\Lambda_G$ is defined and it is cobounded by the above.

Since $\Mod(S)$ possesses a torsion free subgroup of finite index \cite{serretf}, so does $G$, and so, to complete the proof, it suffices to show that $G$ contains no reducible element. But iterating such an element on $\Lambda_G$ would produce a non-filling lamination in $\Lambda_G$, which is excluded by the above.
\end{proof}


\section{Hulls}

\subsection{Minsky's quasi-projections to Teichm\"uller geodesics}

Following Minsky \cite{minskycrelle}, given a closed set $\mathrm C$  in $\T(S)$, we define a closest--points projection from $\T(S)$ to the set of subsets of $\mathrm C$,
\[
\pi_{\mathrm C} \co \T(S) \to \mathcal P (\mathrm{C}),
\]
by demanding that
\[
\pi_\mathrm{C}(X)= \{Y \in \mathrm{C}\, |\, \d_\T(X,Y)=\d_\T(X,\mathrm{C})\}
\]
where $\d_\T(X,\mathrm{C})=\inf_{Y \in \mathrm{C}} \d_\T(X,Y)$.  Given a set $\mathcal X \subset \T(S)$, we abuse notation and refer to $\cup_{X\in \mathcal X} \pi_\mathrm{C}(X) \subset \mathrm{C}$ as $\pi_\mathrm{C}(\mathcal X)$.

Minsky has proven that these projections behave in much the same way as such projections in $\mathbb H^3$.

\begin{theorem}[Contraction Theorem \cite{minskycrelle}]\label{minskycontraction} For every $\epsilon > 0$ there is a constant $b$, depending only on $\epsilon$ and the topological type of $S$ such that for any $\epsilon$--cobounded geodesic $\tau$ and $X$ in $\T(S)$, \[ \mathrm{diam}(\pi_\tau(\N_{\mathrm{d}_\T(X,\tau)}(X))) \leq b. \] \end{theorem}

\begin{theorem}[Corollary 4.1 of \cite{minskycrelle}]\label{quasiproj}  For every $\epsilon > 0$ there is a constant $b$, depending only on $\epsilon$ and the topological type of $S$ such that the following holds. Let $\tau$ be an $\epsilon$--cobounded geodesic in $\T(S)$. If $R > 0$ and points $X,Y \in \T(S)$ are connected by a path of length $T$ that remains outside an $R$--neighborhood of $\tau$, then 
\[
\mathrm{diam}(\pi_\tau(X) \cup \pi_\tau(Y)) \leq \frac{b}{R} T + b.
\]
Furthermore, for any $X,Y \in \T(S)$,
\[
\mathrm{diam}(\pi_\tau(X) \cup \pi_\tau(Y)) \leq \d_\T(X,Y) + 4b.
\]
\end{theorem}

\begin{theorem}[Theorem 4.2 of \cite{minskycrelle}]\label{quasistab} For every $\epsilon > 0$, $K  \geq 1$ and $C \geq 0$, there exists a constant $D$ depending on $\epsilon, K, C$ and the topological type of $S$ such that the following holds.  Let $\psi$ be a $(K, C)$--quasi-geodesic path in $\T(S)$ whose endpoints in $\T(S)$ are connected by an $\epsilon$--cobounded geodesic $\tau$.  Then $\psi$ remains in the $D$--neighborhood of $\tau$. 
\end{theorem}

\subsection{Thin triangles and the hull}\label{trianglesect}

We need the following general fact about cobounded geodesic triangles in Teichm\"uller space.
A different proof has been discovered by M. Duchin, see \cite{duchin}.

\begin{theorem}[Thick triangles are thin]\label{thickthin} For every $\epsilon > 0$ there is a $\delta > 0$ such that if $\triangle$ is a geodesic triangle with vertices in $\T(S) \cup \PML(S)$ whose sides are $\epsilon$--cobounded, then $\triangle$ is $\delta$--thin.
\end{theorem}
\begin{proof} Let $X$, $Y$, and $Z$ be points in $\T(S)$.  Let $P$ be a point in the geodesic segment $[X,Y]$ that minimizes the distance between $Z$ and that segment.  It is shown in the proof of Lemma 7.2 of \cite{MM1} that the path $[X,P] \cup [P,Z]$ is a $(3,0)$--quasi-geodesic---in fact, this is true in any geodesic metric space.  We include the proof here for the reader's convenience.

If $P$ and $X$ coincide, there is nothing to do.  So suppose that $P \neq X$. Let $U$ be a point in $[X,P]$, $V$ a point in $[P,Z]$. By the triangle inequality and choice of $P$, $P$ also minimizes the distance between $V$ and $[X,Y]$.  So $\d_\T(U,V) \geq \d_\T(P,V)$.

By the triangle inequality, $\d_\T(U,V) \geq \d_\T(U,P) - \d_\T(P,V)$.  Together with the previous inequality, we have $3 \d_\T(U,V) \geq \d_\T(U,P) + \d_\T(P,V)$.  So, $[X,P] \cup [P,Z]$ is indeed a $(3,0)$--quasi-geodesic.

Now suppose that the geodesic segments $[X,Y]$, $[X,Z]$, and $[Y,Z]$ are $\epsilon$--cobounded.

By Theorem \ref{quasistab}, the path $[X,P] \cup [P,Z]$ remains in a $D(3,0,\epsilon)$--neighborhood of $[X,Z]$. By symmetry, the path $[Y,P] \cup [P,Z]$ lies in the $D(3,0,\epsilon)$--neighborhood of $[Y,Z]$.  In particular, the segment $[X,Y]$ lies in the $D(3,0,\epsilon)$--neighborhood of the union \linebreak $[X,Z] \cup [Y,Z]$.  Symmetry guarantees that the triangle $\triangle[X,Y,Z]$ is $D(3,0,\epsilon)$--thin.

We continue to assume that the sides of the triangle are cobounded and turn to the case where at least one of $X$, $Y$ and $Z$ lie in $\PML(S)$.  Suppose that $W \in \{X,Y,Z\}$ is such that $W \in \PML(S)$ and let $W'$ and $W''$ be points lying in the interiors of the sides incident to $W$ as pictured in Figure \ref{truncate}. Since the ray $[W',W)$ is cobounded, $W$ is uniquely ergodic, and by Masur's Asymptotic Rays Theorem, the rays $[W',W)$ and $[W'',W)$ are asymptotic. We re-choose $W'$ and $W''$ so that $\d_\T(W',W'') \leq 1$ and the rays $[W',W)$ and $[W'',W)$ are contained in the $1$--neighborhoods of each other.

We truncate the triangle $\triangle[X,Y,Z]$ at any such $W$ to obtain a geodesic polygon inscribed with the triangle $\triangle[X',Y',Z']$, where $W'=W$ if $W \in \T(S)$---the possibilities are depicted in Figure \ref{truncate}.  The result is composed of the triangle $\triangle[X',Y',Z']$ and at most three geodesic triangles of a special type: each has an $\epsilon$--cobounded side, a side of length at most $1$, and a side of $\triangle[X',Y',Z']$.  The union of the latter two sides is a $(1,1)$--quasi-geodesic---as is the union of the former two---and since the remaining side is $\epsilon$--cobounded, this union is contained in the $D(1,1, \epsilon)$--neighborhood of that side, by Theorem \ref{quasistab}. Moreover, this implies that the sides of $\triangle[X',Y',Z']$ are $\epsilon'$--cobounded, for some $\epsilon'$ depending only on $\epsilon$ and the topological type of $S$. 
\begin{figure}
\begin{center}
\input{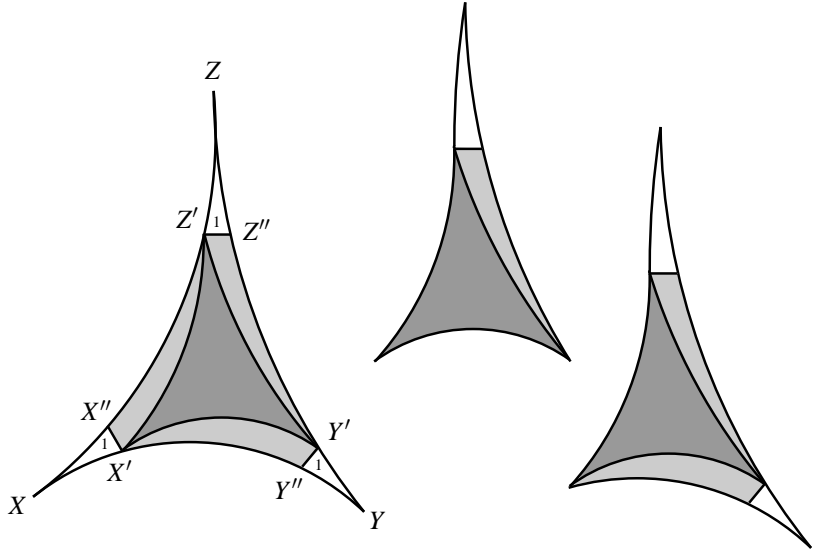}
\end{center}
\caption{Three potential ideal triangles, their polygons, and inscribed triangles.}
\label{truncate}
\end{figure}
By our previous argument, the triangle $\triangle[X',Y',Z']$ is $D(3,0,\epsilon')$--thin. Yet another application of Theorem \ref{quasistab} shows that the special triangles above are $M$--thin, where $M=\max \{D(1,1,\epsilon'), D(1,1,\epsilon)\}$.

In any case, it is easily seen that our triangle $\triangle[X,Y,Z]$ is $\delta$--thin for any $\delta$ larger than $2M + D(3,0,\epsilon') + 1$.
\end{proof}

\begin{theorem}\label{convexhull} For every $\epsilon > 0$, there is an $A \geq 0$ such that if a weak hull $\WH = \WH_\A$ is $\epsilon$--cobounded, then $\WH \cup \A$ is $A$--quasi-convex.  Moreover, any two points in $\WH$ are within $2 \delta$ of a bi-infinite geodesic in $\WH$ (where $\delta$ is given by Theorem \ref{thickthin}).
\end{theorem}
\begin{proof} Note that since $\WH$ is cobounded, every lamination in $\A$ is uniquely ergodic by Masur's criterion and so the end of any geodesic ray in $\WH$ converges in \linebreak $\T(S) \cup \PML(S)$ to its direction, by the Two Boundaries Theorem.  In particular, every pair of distinct points in $\A$ binds $S$.

Let $\delta$ be the constant given by Theorem \ref{thickthin} and let $X$ and $Y$ be two points in $\WH$. 

We begin by finding a bi-infinite geodesic $\gamma$ in $\WH$ such that $X$ and $Y$ are both within $2\delta$ of $\gamma$. If $X$ and $Y$ lie in a geodesic contained in $\WH$, we are done. Otherwise, there are two bi-infinite geodesics $\sigma$ and $\tau$ in $\WH$ containing $X$ and $Y$, respectively. 

There are two cases to consider, when $\sigma$ and $\tau$ have an endpoint in common, and when they do not.

In the first case, $\sigma$ and $\tau$ are two sides of a geodesic triangle contained in $\WH$.  By Theorem \ref{thickthin}, this triangle is $\delta$--thin, and the desired geodesic is easily found.

In the second case, $\sigma$ and $\tau$ determine four points in $\PML(S)$. Join every pair of these points by a Teichm\"uller geodesic.  The resulting union of geodesics in \linebreak $\T(S) \cup \PML(S)$ is the $1$--skeleton of a $3$--simplex, and we refer to the four geodesic triangles in this configuration as the \textbf{faces of the simplex}, the geodesics themselves as \textbf{edges}.  All of the edges are contained in $\WH$ and so all of the faces of the simplex are $\delta$--thin.

Since $\sigma$ and $\tau$ are each incident to two faces of the simplex, for $W \in \{X,Y\}$ there are at least three edges at a distance at most $2\delta$ from $W$.  If for some $W$ there are four edges at such a distance, we know that $X$ and $Y$ are $2\delta$ away from a common edge, by the pigeon--hole principle.  Now, it is easy to see that if for each $W \in \{X,Y\}$ there are exactly three edges a distance at most $2\delta$ from $W$, then these three edges share a vertex. But two such triples of edges in a $3$-simplex must share a common edge.

In any case, the common edge is the desired geodesic $\gamma$, and the second half of the theorem follows.

Joining the geodesic segment $[X,Y]$ to $\gamma$ by geodesic segments yields a $(1,4\delta)$--quasi-geodesic, which must lie in the $D(1,4\delta,\epsilon)$--neighborhood of $\gamma$, where $D(1,4\delta,\epsilon)$ is the constant given by Theorem \ref{quasistab}.  In particular, the segment $[X,Y]$ lies in the $D(1,4\delta,\epsilon)$--neighborhood of $\WH$.

If both $X$ and $Y$ lie in $\PML(S)$, they are the negative and positive directions of a geodesic contained in $\WH$. 

If $Y$, say, lies in $\PML(S)$ and $X$ lies in $\T(S)$, let $\sigma$ be a bi-infinite geodesic in $\WH$ containing $X$.  Joining the endpoints of $\sigma$ to $Y$ by geodesics in $\WH$ yields a triangle that is $\delta$--thin.  So, one of the geodesics containing $Y$ is within $\delta$ of $X$.  Call this geodesic $\gamma$ and let $Z \in \gamma$ be within $\delta$ of $X$. By the Asymptotic Rays Theorem, $[X,Y)$ and $[Z,Y)$ are asymptotic and so we may choose points $X'$ and $Z'$ on these rays, respectively, so that $\d_\T(X',Z') \leq 1$ and the rays $[X',Y)$ and $[Z',Y)$ are contained in the $1$--neighborhoods of each other.  The path $[Z,X] \cup [X,X'] \cup [X',Z']$ is a $(1,\delta + 1)$--quasi-geodesic and so, by Theorem \ref{quasistab},  it lies in the $D(1,\delta + 1,\epsilon)$--neighborhood of $\gamma$.  We conclude that the ray $[X,Y)$ lies in the $(D(1,\delta + 1,\epsilon)+1)$--neighborhood of $\gamma$.

Setting $A = D(1,4\delta,\epsilon) + D(1,\delta + 1,\epsilon)+1$ completes the proof.
\end{proof}

\begin{corollary}\label{compactnessalone} If a subgroup $G$ of $\Mod(S)$ has a limit set whose weak hull $\WH_G$ is defined and $G$ acts cocompactly on $\WH_G$, then $G$ is convex cocompact. 
\end{corollary}
\begin{proof} If $G$ acts cocompactly on $\WH_G$, then a $G$--orbit in $\WH_G$ is $B$--dense for some positive number $B$.  Since $\WH_G$ is $A$--quasi-convex for some $A$, we conclude that the $G$--orbit is $(A+B)$--quasi-convex. 
\end{proof}

We also have

\begin{corollary}\label{hyperbolichull} Let $\WH$ be an $\epsilon$--cobounded weak hull and let $A$ be the constant given by Theorem \ref{convexhull}. Let $\WH^A =\N_A(\WH)$ equipped with the induced path metric. Then $\WH^A$ is a proper $\delta$--hyperbolic metric space for some $\delta$.
\end{corollary}
\begin{proof} By the choice of $A$, the restriction of the metric on $\WH^A$ to $\WH$ agrees with the restriction of the Teichm\"uller metric and every geodesic triangle with vertices in $\WH$ has Teichm\"uller geodesic edges.  Let $\triangle$ be such a triangle. Since $\WH$ is $\epsilon$--cobounded and $\triangle$ is contained in its $A$--neighborhood, $\triangle$ is $\epsilon'$--cobounded for some $\epsilon'$.  By Theorem \ref{thickthin}, there is a $\delta'$ depending only on $\epsilon'$ and $S$ such that $\triangle$ is $\delta'$--thin. As $\WH$ is $A$--dense in $\WH^A$, we conclude that $\WH^A$ is $\delta$--hyperbolic for some $\delta$.
\end{proof}

\section{Kleinian manifolds} \label{kleincccsect}

The following is part of Proposition 5.1 of \cite{klarreich}.

\begin{proposition}[Klarreich]\label{klar5.1} Let $X_n$ and $Y_n$ be sequences in $\T(S)$ that converge in Teichm\"uller's compactification to filling laminations $[\mu]$ and $[\lambda]$. Let $\tau_n$ be the sequence of Teichm\"uller geodesic segments joining $X_n$ to $Y_n$.  Then for every accumulation point $[\nu]$ of $\{\tau_n\}$ in Teichm\"uller's boundary one has $i(\nu,\mu) = 0$ or $i(\nu,\lambda) = 0$.\qed
\end{proposition}

We need

\begin{proposition}\label{klarprop1} Let $\mathfrak A \subset \PML(S)$ be a closed set consisting entirely of uniquely ergodic laminations. Let $[\mu_n]$ and $[\lambda_n]$ be sequences in $\A$ converging to $[\mu]$ and $[\lambda]$, and, for each $n$, let $\tau_n$ be the Teichm\"uller geodesic with negative and positive directions $[\mu_n]$ and $[\lambda_n]$.  Then the set of accumulation points  of $\{\tau_n\}$ in Thurston's boundary is contained in $\{[\mu],[\lambda] \}$. 
\end{proposition}
\begin{proof} By Masur's Two Boundaries Theorem, the ends of the geodesics $\tau_n$ converge to their directions $[\mu_n]$ and $[\lambda_n]$.  The proposition now follows from Proposition \ref{klar5.1} by a diagonal argument and the unique ergodicity of $[\mu]$ and $[\lambda]$. 
\end{proof}

\begin{lemma}\label{thehullisclosed} Let $G$ be a subgroup of $\Mod(S)$ such that every lamination in $\Lambda_G$ is uniquely ergodic.  Then the weak hull $\WH_G$ of $\Lambda_G$ is closed in $\T(S) \cup \Delta_G$.  In particular, $\WH_G \cup \Lambda_G$ is closed in $\T(S) \cup \PML(S)$.
\end{lemma}
\begin{proof} Let $X_n$ be a sequence in $\WH_G$ and let $\tau_n$ be a sequence of bi-infinite geodesics in $\WH_G$ containing the $X_n$.  

Suppose that the $X_n$ accumulate at a point $X$ in $\T(S)$. We may pass to a subsequence so that the $X_n$ converge to $X$.  The Arzel\`a--Ascoli Theorem allows us to pass to a further subsequence so that the $\tau_n$ converge to a geodesic through $X$. Since $\Lambda_G$ is closed, the limiting geodesic lies in the weak hull $\WH_G$.

If the $X_n$ accumulate at a point $[\nu]$ in $\PML(S)$, pass to a subsequence so that the $X_n$ converge to $[\nu]$ and so that the ends of the geodesics $\tau_n$ converge to projective measured laminations $[\lambda]$ and $[\mu]$. By Proposition \ref{klarprop1}, $[\nu] \in \{[\mu],[\lambda]\}\subset \Lambda_G$. 
\end{proof}

\begin{theorem}\label{kleinccc} If $G$ is a subgroup of $\Mod(S)$ that acts cocompactly on $\T(S) \cup \Delta_G$, then $G$ is convex cocompact. 
\end{theorem}
\begin{proof} Suppose that $\dot M_G = (\T(S) \cup \Delta_G)/G$ is compact. Then $\Delta_G/G$ is compact. By Theorem \ref{domainsub}, every lamination in $\Lambda_G$ is uniquely ergodic, $Z\Lambda_G = \Lambda_G$, and the weak hull $\WH_G$ of the limit set $\Lambda_G$ is defined.  By Lemma \ref{thehullisclosed}, $\WH_G$ is closed in $\T(S) \cup \Delta_G$.

As $G$ acts cocompactly and properly discontinuously on $\T(S) \cup \Delta_G$ and $\WH_G$ is closed therein, $G$ acts cocompactly on $\WH_G$.  The theorem now follows from Corollary \ref{compactnessalone}.
\end{proof}

\section{Hulls revisited: quasi-projections}\label{hullprojsect}

With suitably adjusted constants, Minsky's quasi-projection theorems hold for cobound-ed weak hulls.

\begin{theorem}[Hull contraction]\label{contraction} Given $\epsilon >0$ there is a constant $c$ depending only on $\epsilon$ and the topological type of $S$ such that for any $\epsilon$--cobounded weak hull $\WH$ and any point $X$ in $\T(S)$,
\[
\mathrm{diam}(\pi_\WH(\N_L(X))) \leq c,
\]
where $L = \d_\T(X,\WH)$.
\end{theorem}
\begin{proof} Fix an $\epsilon$--cobounded weak hull $\WH$ and let $X$ be a point in $\T(S)$. By Theorem \ref{convexhull}, $\WH$ is $A$--quasi-convex for some $A$.

If $L  = \d_\T(X,\WH)< A$, $\mathrm{diam}(\pi_\WH(\N_L(X)))$ is at most $2A$, and so we suppose that $L \geq A$.

We begin by bounding the distance between two points in $\pi_\WH(X)$.  Strictly speaking, this follows from the proof of the theorem given below. As it is a basic ingredient in the proof,  we include the argument in the interest of clarity.  We refer the reader to Figure \ref{projfig} for a diagram of the following.
\begin{figure}
\begin{center}
\input{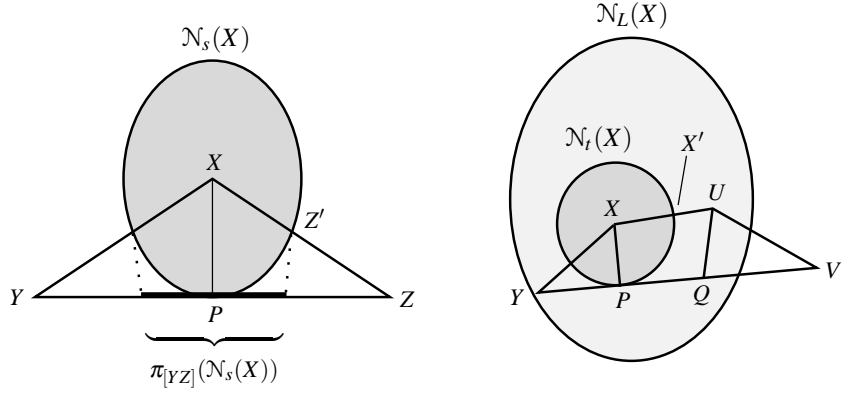}
\end{center}
\caption{Bounding the diameter of the projection.}
\label{projfig}
\end{figure}

Let $Y$ and $Z$ be two points in $\pi_\WH(X)$. The geodesic segment $[Y,Z]$ is contained in the $A$--neighborhood of $\WH$. Let $P$ be a point on $[Y,Z]$ minimizing the distance between $X$ and that segment, and let $s=\d_\T(X,P)$. Note that $s \geq L-A$, for if not, $\d_\T(X,\WH) < L$ as $[Y,Z] \subset \N_A(\WH)$ and this contradicts the fact that $\d_\T(X,\WH) = L$.

By Minsky's Contraction Theorem, there is a constant $b$ depending only on $\epsilon$ and the topological type of $S$ such that
\[
\mathrm{diam}(\pi_{[Y,Z]}(\N_s(X))) \leq b.
\]

Consider the geodesic segment $[X,Z]$. This is composed of a segment $[X,Z']$ of length $s$ and a segment $[Z',Z]$ of length at most $A$, as the whole segment has length $L$. By the above, the diameter of $\pi_{[Y,Z]}([X,Z'])$ is at most $b$.  By Theorem \ref{quasiproj}, the diameter of $\pi_{[Y,Z]}([Z',Z])$ is at most $A + 4b$. So the diameter of $\pi_{[Y,Z]}([X,Z])$ is at most $A + 5b$.

Since $P \in \pi_{[Y,Z]}(X)$ and $\pi_{[Y,Z]}(Z) = \{Z\}$, we conclude that $\d_\T(P,Z) \leq A + 5b$. Symmetry yields the same bound for the distance between $P$ and $Y$, and we conclude that $\d_\T(Y,Z) \leq 2A + 10b$. 

We now turn to the proof of the theorem, continuing to let $Y$ denote a point in $\pi_\WH(X)$, letting $U$ be a point different from $X$ in $\N_L(X)$, and letting $V$ be an element of $\pi_\WH(U)$.

Let $Q$ be a point of $[Y,V]$ minimizing the distance between $U$ and that segment.  An argument similar to the above shows that $\d_\T(Q,V) \leq A + 5b$.  Let $P$ be a point on $[Y,Q]$ minimizing the distance between $X$ and $[Y,Q]$.  As $[Y,V] \subset \N_A(\WH)$, so is $[Y,Q]$, and we have $t = \d_\T(X,P) \geq L-A$.  Again, we have that $\d_\T(Y,P) \leq A+5b$.

Now, since $U \in \N_L(X)$, $\d_\T(X,U) \leq L$ and so $[X,U]$ is composed of two segments $[X,X']$ and $[X',U]$, where $[X,X']$ lies in the $(L-A)$--neighborhood of $X$ and $[X',U]$ has length at most $A$.  We conclude that 
\[
\mathrm{diam}(\pi_{[P,Q]}([X,X'])) \leq b,
\]
by Minsky's Contraction Theorem, and that
\[
\mathrm{diam}(\pi_{[P,Q]}([X',U])) \leq A + 4b,
\]
by Theorem \ref{quasiproj}.  So, the distance between $Y$ and $V$ is at most $3A + 15b$, and we conclude that 
\[
\mathrm{diam}(\pi_\WH(\N_L(X))) \leq 6A + 30b.
\]

Letting $c = 6A + 30b$ completes the proof.
\end{proof}

As in \cite{minskycrelle}, this easily yields analogs of Theorems \ref{quasiproj} and \ref{quasistab} for weak hulls, which we state for completeness.

\begin{theorem}  For every $\epsilon > 0$ there is a constant $c$, depending only on $\epsilon$ and the topological type of $S$ such that the following holds. Let $\WH$ be an $\epsilon$--cobounded weak hull in $\T(S)$. If $R > 0$ and points $X,Y \in \T(S)$ are connected by a path of length $T$ that remains outside an $R$--neighborhood of $\WH$, then 
\[
\mathrm{diam}(\pi_\WH(X) \cup \pi_\WH(Y)) \leq \frac{c}{R} T + c.
\]
Furthermore, for any $X,Y \in \T(S)$,
\[
\mathrm{diam}(\pi_\WH(X) \cup \pi_\WH(Y)) \leq \d_\T(X,Y) + 4c.
\]
\end{theorem}

\begin{theorem}\label{hullquasistab} Let $\psi$ be a $(K, C)$--quasi-geodesic path in $\T(S)$ whose endpoints in $\T(S)$ are contained in an $\epsilon$--cobounded weak hull $\WH$.  Then $\psi$ remains in a $E(K, C, \epsilon)$ neighborhood of $\WH$. 
\end{theorem}

\begin{remark} Theorem \ref{hullquasistab} also follows directly from Theorem \ref{quasistab} and Theorem \ref{convexhull}.
\end{remark}

\section{Quasi-isometric embedding in $\C(S)$} \label{qisect}

Let $G$ be a subgroup of $\Mod(S)$ with finite generating set $\mathcal U$ and word metric $\d_{\mathcal U}$.
For any $v \in \C(S)$, the $G$-orbit $G v$ of $v$ defines a map $\Phi_v\co G \to \C(S)$.
We have the following

\begin{theorem} \label{ifqithencc}
For any $v \in \C(S)$, if $\Phi_v$ is a quasi-isometric embedding, then $G$ is convex cocompact.
\end{theorem}
\begin{remark}
See also U. Hamenst\"adt \cite{hamenstadt}.
\end{remark}

If $\Phi_v$ is a $(K,C)$--quasi-isometric embedding, then for any $u$ in $\C(S)$, the map $\Phi_u$ is a $(K, C')$--quasi-isometric embedding, where $C' = C + 2\d_\C(u,v)$---in particular, we may assume that $v$ is any point of $\C(S)$, when a choice of $v$ is convenient.

Given a point $X$ in $\T(S)$, the $G$--orbit $G X$ of $G$ defines a map $\Psi_X\co G \to \T(S)$.  It so happens that $\Phi_v$ being a quasi-isometric embedding implies that $\Psi_X$ is as well.  We record this in the following


\begin{lemma} \label{Kaleb} If $\Phi_v$ is a quasi-isometric embedding for some $v \in \C(S)$, then for any point $X$ in $\T(S)$, the map $\Psi_X\co G \to \T(S)$ is a quasi-isometric embedding.
\end{lemma}
\begin{proof}  Since $\Psi_X$ is defined by taking an element $h$ to $h X$, the desired upper bound is an immediate consequence of the finite generation of $G$.

The Teichm\"uller space sits naturally in the electric Teichm\"uller space $\Tel(S)$, see the proof of Theorem \ref{ifccthenqi} for a definition.
By Lemma 3.1 of \cite{MM1}, and its proof, the electric space and $\C(S)$ are $\Mod(S)$--equivariantly quasi-isometric.
We may assume that $v \in \C(S)$ is the image of $X$ under such a quasi-isometry.
Since $\Phi_v$ is a quasi-isometric embedding and the inclusion $\T(S) \to \Tel(S)$ is Lipschitz, we obtain the desired lower bound.
\end{proof}

\subsection{The boundary and ending laminations}

\noindent By Theorem 1.1 of \cite{MM1}, see also \cite{bowditch}, $\C(S)$ is $\delta'$--hyperbolic for some $\delta'$.
If \linebreak $\Phi_v \co G \to \C(S)$ is a quasi-isometric embedding then $G$ is $\delta$--hyperbolic for some $\delta$, the map $\Phi_v$ has a continuous extension
\[
\Phi_v \co G \cup \partial G \to \C(S) \cup \partial \C(S),
\]
and the restriction
\[
\partial \Phi_v \co \partial G \to \partial \C(S)
\]
is a topological embedding, see Th\'eor\`eme 2.2 of \cite{coor}.

By the stability of quasi-geodesics in $\delta$--hyperbolic metric spaces, see Th\'eor\`eme 1.2 of \cite{coor}, there exists an $A > 0$ such that for any geodesic $\g$ in $G$, the quasi-geodesic $\Phi_v(\g)$ and any geodesic joining its endpoints have Hausdorff distance at most $A$.
In particular, for any distinct pair of points $x,y \in \Phi_v(G \cup \partial G)$, any geodesic between $x$ and $y$ is contained in $\N_{A}(\Phi_v(G))$---thus $\Phi_v(G)$ is $A$-quasi-convex.

In the next section, we find estimates required to cobound the weak hull (see Corollary \ref{cobound}).
To do this, we must recall the geometric description of $\partial \C(S)$.

Let $\L_{\mathrm{fill}}(S)$ be the set of filling laminations in $\PML(S)$ and let \linebreak $\mathrm F \co \L_{\mathrm{fill}}(S) \to \UML(S)$ be the map that forgets transverse measures.
The image of $\mathrm F$ is the space of potential ending laminations for hyperbolic $3$--manifolds homeomorphic to $S \times \mathbb R$ and is denoted here by $\EL(S)$.
It is a theorem of E. Klarreich  \cite{klarreich} that $\partial \C(S)$ is naturally homeomorphic to $\EL(S)$ so that if a quasi-geodesic limits to $|\mu| \in \partial \C(S)$, then every accumulation point in $\PML(S)$ of its vertices---viewed as elements in $\PML(S)$---projects to $|\mu|$ under $\mathrm F$.
In particular, for any $m \in \partial G$, $\Phi_v(m)$ is naturally identified with a lamination in $\EL(S)$.

\subsection{Bounding the subsurface projection coefficients}

We make repeated use of the following theorem of Masur and Minsky \cite{MM2}.

\begin{theorem}[Masur--Minsky] \label{mmbgi}
There exists a constant $M = M(S)$ with the following property.
Let $Y$ be a proper domain of $S$ with $\xi(Y) \neq 3$ and let $\gamma$ be a geodesic segment, ray, or bi-infinite line in $\C(S)$, such that $\pi_{Y}(v) \neq \emptyset$ for every vertex $v$ of $\gamma$.
Then
\[
\diam_{Y}(\gamma) \leq M.
\]
\end{theorem}

The main theorem allowing us to cobound the hull is 

\begin{theorem}[Quasi-isometric projection bound]\label{proj} If $\Phi_v$ is a quasi-isometric embedding for some $v \in \C(S)$, then there exists a constant $D > 0$ such that for any two distinct points $m_-,m_+ \in \partial G$ and any proper domain $Y \subset S$ with $\xi(Y) \neq 3$, we have
\[
\d_Y(\Phi_v(m_-), \Phi_v(m_+)) \leq D.
\]
\end{theorem}
\noindent As the proof is technical, we pause to sketch the argument.

Given distinct points $m_-,m_+$ in $\partial G$, there is a geodesic in $G$ joining them. This geodesic is carried to a quasi-geodesic in $\C(S)$, which is uniformly close to a geodesic $\gamma$ joining $\Phi_v(m_-)$ and $\Phi_v(m_+)$. To bound a coefficient $\d_Y(\Phi_v(m_-), \Phi_v(m_+))$, it suffices to bound $\diam_{Y}(\gamma)$. 

If $Y$ is a proper domain whose boundary is far from $\gamma$, Theorem \ref{mmbgi} provides a bound on $\diam_{Y}(\gamma)$.  If $\partial Y$ is close to $\gamma$, it is close to $\Phi_v(G)$. In fact, we may assume that $\partial Y$ is close to $\Phi_v(\1)$ by translating. Since the two ends of $\gamma$ diverge,  $\gamma$ may be decomposed into three parts: a finite segment $\gamma_0$ near $\partial Y$ and two infinite rays $\gamma_-$ and $\gamma_+$ far from $\partial Y$.  Theorem \ref{mmbgi} again bounds $\diam_{Y}(\gamma_\pm)$. The segment $\gamma_0$ fellow travels the image of a geodesic segment in $G$ lying in a fixed neighborhood of $\1$. Finiteness of this neighborhood allows us to bound $\diam_Y(\gamma_0)$.  The triangle inequality provides the bound on $\diam_Y(\gamma)$.
\begin{proof} 
We let $\Omega$ denote the set of pairs of distinct points in $\partial G$:
\[
\Omega = \{(m_-,m_+) \, | \, m_-, m_+ \in \partial G \mbox{ and } m_- \neq m_+ \}.
\]

We assume $\Phi_v$ is a $(K,C)$--quasi-isometry and as noted above, $\Phi_v(G)$ is $A$--quasi-convex.
It is convenient to assume further that we have chosen $A$ sufficiently large so that for any geodesic $\g$ in $G \cup \partial G$ and any geodesic $\gamma$ connecting the endpoints of $\Phi_v(\g)$, any closest point projection map from $\Phi_v(\g)$ to $\gamma$ is \textbf{$A$--coarsely order preserving}: if $u_0, u_1, u_2 \in \gamma$ are pairwise separated by a distance at least $A$ and $u_0 < u_1 < u_2$, then for every triple $h_0,h_1,h_2 \in \g$ for which $\Phi_v(h_i)$ is a point closest to $u_i$, $i = 0,1,2$, we have $h_0 < h_1 < h_2$.

We partition the proper domains of $S$ into two classes
\begin{align*}
\Dom(\mathrm{far})  =  \{ Y \subset S \, | \, \d_\C(\partial Y,\Phi_v(G)) \geq A + 2 \},\\
\Dom(\mathrm{near})  =  \{ Y \subset S \, | \, \d_\C(\partial Y,\Phi_v(G)) < A + 2\}, 
\end{align*}
and define
\[
\Dom(0)  =  \{ Y \subset S \, | \, \d_\C(\partial Y,v) < A + 2\} \subset \Dom(\mathrm{near}).
\]

Let $Y$ be an element of $\Dom(\mathrm{far})$. By Lemma 5.14 of \cite{ELCI}, for any pair $(m_-,m_+)$ in $\Omega$, there exists a geodesic $\gamma$ between $\Phi_v(m_-)$ and $\Phi_v(m_+)$.  By our choice of $A$, this lies in ${\mathcal N}_{A}(\Phi_v(G))$ and so $\d_\C(\partial Y,\gamma) \geq 2$. In particular, $\pi_Y(u) \neq \emptyset$ for every $u \in \gamma$, and Theorem \ref{mmbgi} implies that
\[
\mathrm{diam}_Y(\gamma) \leq M.
\]
Let $\{v_i\}$ be the vertices of $\gamma$. Since $\Phi_v(m_+)$ is a sub-lamination of any accumulation point of $\{v_i \}_{i=0}^{\infty}$ in the Hausdorff topology, it follows that for any subsurface $Y$, we have
\[
\pi_Y(\Phi_v(m_+)) \subset \pi_Y(v_i)
\]          
for sufficiently large $i$. Along with a similar argument for $\Phi_v(m_-)$, this yields
\[
\diam_Y(\gamma) = \diam_Y(\{v_i\} \cup \{\Phi_v(m_\pm)\}),
\]
and we conclude that
\[
\d_Y(\Phi_v(m_-),\Phi_v(m_+)) \leq M.
\]

Now let $Y$ be a domain in $\Dom(\mathrm{near})$, $(m_-,m_+)$ a pair in $\Omega$, and $\g$ a geodesic in $G$ joining $m_-$ to $m_+$. Let $h \in \g$ be a point for which $\d_\C(\partial Y,\Phi_v(h)) < A + 2$.
Then $h^{-1}(\partial Y) = \partial (h^{-1}(Y))$, $\1 \in h^{-1}(\g)$, and $h^{-1}(Y) \in \Dom(0)$.
Furthermore, $h^{-1}(\g)$ has endpoints $h^{-1}(m_-)$ and $h^{-1}(m_+)$ and 
\[
\d_Y(\Phi_v(m_-),\Phi_v(m_+)) = \d_{h^{-1}(Y)}(\Phi_v(h^{-1}(m_-)),\Phi_v(h^{-1}(m_+))).
\]
So it suffices to find a constant $D'$ such that
\begin{equation} \label{whatweneed}
\d_Y(\Phi_v(m_-),\Phi_v(m_+)) \leq D'
\end{equation}
whenever $Y \in \Dom(0)$ and $(m_-,m_+)$ is a pair joined by a geodesic through $\1$. Setting $D = \max\{ D',M\}$ will complete the proof.\\

\noindent \textbf{Finding $D'$.} We fix a constant $R$ satisfying
\begin{equation} \label{reqn}
R \geq K (4 A + 5 + C)
\end{equation}
and refer the reader to Figure \ref{quasigeo1} for a schematic of what follows.

\begin{figure}
\input{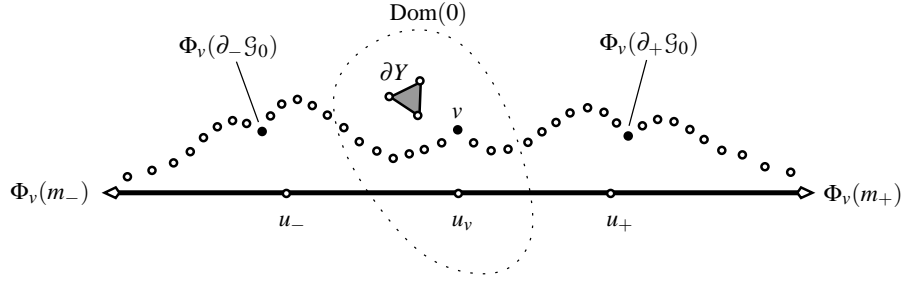}
\caption{The quasi-geodesic $\Phi_v(\g)$, its fellow--traveling geodesic $\gamma$, and some points of interest.}
\label{quasigeo1}
\end{figure}

We fix a pair $(m_-,m_+)$ in $\Omega$ and a geodesic $\g$ through $\1$ joining $m_-$ and $m_+$.
Let $\g_0 \subset \g$ denote the intersection of $\g$ with the ball of radius $R$ about $\1$.

Next, let $\partial_-\g_0$ and $\partial_+\g_0$ denote the initial and terminal points of $\g_0$, respectively.
Since $\Phi_v(\1) = v$ we have
\begin{equation}\label{eqn1}
\begin{split}
\d_\C(\Phi_v(\partial_\pm \g_0),v) & \geq \frac{1}{K} \d_U(\partial_\pm \g_0,\1) - C \\
& = \frac{R}{K} - C \\
& \geq  4 A + 5.
\end{split}
\end{equation}
Similarly, we observe
\begin{equation} \label{eqn2}
\begin{split}
\d_\C(\Phi_v(\partial_-\g_0),\Phi_v(\partial_+\g_0)) & \geq 8 A + 10 + C \\
& \geq 8 A + 10.
\end{split}
\end{equation}
Again by \cite{ELCI}, there is a geodesic $\gamma$ with endpoints $\Phi_v(m_-)$ and $\Phi_v(m_+)$.
This has Hausdorff distance at most $A$ from $\Phi_v(\g)$, and we let $u_-=u_-(\gamma)$ and $u_+=u_+(\gamma)$ denote a pair of closest points on $\gamma$ to $\Phi_v(\partial_-\g_0)$ and $\Phi_v(\partial_+\g_0)$, respectively.
Note that
\begin{equation} \label{eqn3}
\d_\C(u_\pm,\Phi_v(\partial_\pm \g_0)) \leq A. 
\end{equation}
By (\ref{eqn1}) and (\ref{eqn3}) (and the triangle inequality) we have 
\begin{equation} \label{eqn4}
\begin{split}
\d_\C(u_\pm,v) & \geq \d_\C(\Phi_v(\partial_\pm \g_0),v) - \d_\C(u_\pm,\Phi_v(\partial_\pm \g_0))\\
 & \geq 4 A + 5 - A \\
 & = 3 A  + 5.
\end{split} 
\end{equation}
Then, by (\ref{eqn2}) and (\ref{eqn3}),
\begin{equation} \label{eqn5} 
\begin{split}
\d_\C(u_-,u_+) & \geq \d_\C(\Phi_v(\partial_-\g_0),\Phi_v(\partial_+\g_0)) \\
& \quad - \d_\C(\Phi_v(\partial_-\g_0),u_-)\\
& \quad - \d_\C(\Phi_v(\partial_+\g_0),u_+)\\
 & \geq  8 A + 10 - 2 A \\
& = 6 A + 10.
\end{split} 
\end{equation}
In particular, $u_-,u_+$ decomposes $\gamma$ into a pair of geodesic rays $\gamma_{\pm}$ and a geodesic segment $\gamma_0$.
The endpoints of $\gamma_-$, $\gamma_0$, and $\gamma_+$ are $\{\Phi_v(m_-),u_-\}$, $\{u_-,u_+\}$, and \linebreak $\{u_+,$ $\Phi_v(m_+)\}$, respectively.

Let $u_v=u_v(\gamma)$ denote a closest point on $\gamma$ to $v$, which is a distance at most $A$ from $v$.
By (\ref{eqn4}) we have
\begin{equation} \label{eqn6}
\begin{split}
\d_\C(u_\pm,u_v) & \geq \d_\C(u_\pm,v) - \d_\C(u_v,v)\\
&  \geq 3 A + 5 - A \\
& = 2 A + 5. 
\end{split}
\end{equation}
Thus, because any closest point projection to $\gamma$ is $A$--coarsely order preserving, and since $\partial_-\g_0 < \1 < \partial_+\g_0$, it must be that $u_- < u_v < u_+$, and $u_v \in \gamma_0$.

Moreover, by (\ref{eqn6}) and because $\gamma$ is a geodesic, we have, for every $u \in \gamma_{\pm}$,
\begin{equation} \label{eqn7}
\begin{split}
\d_\C(u,v) & \geq  \d_\C(u,u_v) - \d_\C(v,u_v)\\
 & \geq  \d_\C(u_\pm,u_v) - \d_\C(v,u_v)\\
& \geq 2 A + 5 - A \\
&\geq A + 5.
\end{split} 
\end{equation}

Now suppose that $Y \in \Dom(0)$ and $u \in \gamma_{\pm}$.
By (\ref{eqn7}) we have
\begin{equation} \label{eqn8}
\begin{split}
\d_\C(\partial Y,u) & \geq \d_\C(u,v) - (\d_\C(\partial Y,v) + \diam(\partial Y))\\
 & \geq A + 5 - (A + 2 + 1) \\
& = 2
\end{split}
\end{equation}
and therefore $\pi_Y(u) \neq \emptyset$ for every $u \in \gamma_{\pm}$ and hence
\begin{equation} \label{eqn9}
\diam_{Y}(\gamma_{\pm}) \leq M.
\end{equation}
As before, we have
\[
\d_Y(\Phi_v(m_\pm), u_\pm) \leq \mathrm{diam}_Y(\gamma_\pm).
\]

Next, suppose that $\zeta_{\pm}$ are geodesics connecting $\Phi_v(\partial_\pm \g_0)$ to $u_\pm$.
These geodesics have length at most $A$ and therefore for every $u \in \zeta_{\pm}$, by (\ref{eqn1}), we have
\begin{equation} \label{eqn10}
\begin{split}
\d_\C(\partial Y,u) & \geq  \d_\C(v,\Phi_v(\partial_\pm \g_0)) -  \d_\C(u,\Phi_v(\partial_\pm \g_0))\\
 & \quad  - \d_\C(v,\partial Y) - \diam(\partial Y)\\
 & \geq 4 A + 5 - A - (A + 2) -1 \\
&= 2 A + 2. 
 \end{split}
\end{equation}

In particular, we see that $\pi_Y(\zeta_{\pm}) \neq \emptyset$ and so, by Theorem \ref{mmbgi},
\begin{equation} \label{eqn11}
\diam_{Y}(\zeta_{\pm}) \leq M.
\end{equation}

Thus, by (\ref{eqn9}) and (\ref{eqn11}) we obtain
\begin{align*}
\d_Y(\Phi_v(m_-),\Phi_v(m_+)) & \leq  \diam_Y(\gamma_-)  + \diam_Y(\zeta_-)  \\
 & \quad  + \d_Y(\Phi_v(\partial_-\g_0),\Phi_v(\partial_+\g_0))\\
 & \quad + \diam_Y(\zeta_+) + \diam_Y(\gamma_+) \\
 & \leq  4 M + \d_Y(\Phi_v(\partial_-\g_0),\Phi_v(\partial_+\g_0)). 
\end{align*}

Note that this last expression depends only on $Y$, $v$, and $\g_0 \subset \N_R(\1)$.

Since $\Phi_v(\N_R(\1))$ is finite, there is a constant $D''$ such that for each pair $u$ and $w$ in $\Phi_v(\N_R(\1))$, the intersection number $i(u,w)$ is at most $D''$.  As a consequence, there is a constant $D'''$ such that $\d_Y(u,w) \leq D'''$ for all proper domains $Y$.

Setting $D' = 4M + D'''$, condition (\ref{whatweneed}) is satisfied, and the proof is complete.
\end{proof}

\subsection{Proof of Theorem \ref{ifqithencc}}

In \cite{klarreich}, Klarreich shows that $\mathrm F\co \L_{\mathrm{fill}}(S) \to \EL(S)$ is a closed map.
Combining this with Theorem \ref{proj}, we obtain

\begin{corollary}\label{cobound} If $\Phi_v$ is a quasi-isometric embedding, then for any $m \in \partial G$, the lamination $\Phi_v(m)$ is uniquely ergodic.
Moreover, we can lift $\partial \Phi_v$ to a continuous $G$--equivariant map
\[
\partial \Psi\co \partial G \to \Lambda_G \subset \L_{\mathrm{fill}}(S) \subset \PML(S)
\]
parameterizing the limit set, and the weak hull $\WH_G$ of $\Lambda_G$ is defined and cobounded.
\end{corollary}
\begin{proof} Let $m_-$ and $m_+$ be any pair of points in $\partial G$.
If $[\mu_\pm] \in \mathrm F^{-1}(\Phi_v(m_\pm))$, then $\mu_-$ and $\mu_+$ bind $S$ and so correspond to a Teichm\"uller geodesic $\tau(\mu_-, \mu_+)$.
Since $[\mu_\pm]$ have supporting laminations $\Phi_v(m_\pm)$,  Theorem \ref{proj} produces a constant $D$ such that for any proper domain $Y \subset S$ with $\xi(Y) \neq 3$, the subsurface projection coefficient $\d_Y(\mu_-, \mu_+)$ is at most $D$.
By Rafi's Theorem, Theorem \ref{rafitheorem} here, there is an $\epsilon > 0$ depending only on $S$ and $D$ such that the Teichm\"uller geodesic $\tau(\mu_-, \mu_+)$ is $\epsilon$--cobounded.
By Masur's criterion, the laminations $[\mu_\pm]$ (or equivalently, $\Phi_v(m_\pm)$) are uniquely ergodic.
In particular note that $\mathrm F^{-1}(\Phi_v(m_\pm))$ are in fact singletons $\{[\mu_\pm]\}$.
Therefore, for any $m \in \partial G$, the lift $\partial \Psi(m)$ given by $\mathrm F^{-1}(\Phi_v(m))$ is defined.

Since $\partial \Psi$ is the unique lift of $\partial \Phi_v$, it follows that if $V$ is any closed set in $\PML(S)$, $\partial \Phi_v^{-1}(\mathrm{F}(V)) = \partial \Psi^{-1}(V)$.
Because $\mathrm{F}$ is a closed map and $\partial \Phi_v$ is continuous, we see that $\partial \Psi$ is continuous.
Furthermore, as fixed points of hyperbolic elements of $G$ are dense in $\partial G$, fixed points of pseudo-Anosov elements of $G$ are dense in $\partial \Psi(\partial G)$.
Therefore, being the image of a compact set, $\partial \Psi(\partial G)$ must agree with $\Lambda_G$.
In particular, $\WH_G$ is defined and $\epsilon$--cobounded.
\end{proof}

We now give the

\begin{proof}[Proof of Theorem \ref{ifqithencc}] According to Corollary \ref{cobound}, $\WH_G$ is defined and $\epsilon$--cobounded for some $\epsilon$. Let $A$ be the constant given by Theorem \ref{convexhull}.

By Corollary \ref{hyperbolichull}, $\N_A(\WH_G)$ (with the induced path metric) is a $\delta$-hyperbolic metric space for some $\delta$.
Moreover, the inclusion $\WH_G \to \T(S)$ is an isometric embedding
and every geodesic in $\T(S)$ connecting a pair of points in $\WH_G$ is contained in $\N_A(\WH_G)$, by Theorem \ref{convexhull}.

Let $X$ be a point in $\WH_G$. 
By Lemma \ref{Kaleb}, $\Psi_X\co G \to \T(S)$ is a quasi-isometric embedding and so $\Psi_X \co G \to \N_A(\WH_G)$ is a quasi-isometric embedding by the above.
By the stability of quasi-geodesics in a $\delta$--hyperbolic metric space, $\Psi_x(G) \subset \WH_G$ is quasi-convex in $\N_A(\WH_G)$, and so in $\T(S)$.
\end{proof}

\subsection{Electricity and the converse}\label{converse}

\begin{theorem} \label{ifccthenqi}
If $G < \Mod(S)$ is convex cocompact, then for any $v \in \C(S)$, $\Phi_v$ is a quasi-isometric embedding.
\end{theorem}
\begin{remark} See also Hamenst\"adt, \cite{hamenstadt}.
\end{remark}
\begin{proof} Let $\epsilon > 0$. For $\alpha \in \C_0(S)$, let 
\[
\mathrm{thin}(\alpha, \epsilon) = \{ X \in \T(S)\, | \, \ext_X(\alpha)\leq \epsilon \}
\]
and let 
\[
\mathrm{thin}(\epsilon) = \bigcup_{\alpha \in \C_0} \! \mathrm{thin}(\alpha, \epsilon).
\]
Let $G$ be a convex cocompact subgroup of $\Mod(S)$.
Let $\epsilon_0$ be a number small enough so that the nerve of the family $\{\mathrm{thin}(\alpha, \epsilon_0)\}$ is the complex of curves and so that $\WH_G$ is $\epsilon_0$--cobounded.

Recall that $\Tel(S)$ is obtained from $\T(S)$ in the following way.  For each $\mathrm{thin}(\alpha, \epsilon_0)$, we create a new point $u_\alpha$ and adjoin an interval of length $\frac 1 2$ between each point in $\mathrm{thin}(\alpha, \epsilon_0)$ and $u_\alpha$. Taking the induced path metric yields the space $\Tel(S)$.

Since $\C(S)$ is quasi-isometric to $\Tel(S)$, it suffices to show that an orbit in $\Tel(S)$ defines a quasi-isometric embedding $G \to \Tel(S)$. We claim that there are constants $K \geq 1$ and $C \geq 0$ such that the bi-infinite Teichm\"uller geodesics in $\WH_G$ are all $(K, C)$--quasi-geodesics in $\Tel(S)$.
As any two points in $\WH_G$ are within $2\delta = 2\delta(\epsilon_0)$ of a bi-infinite geodesic in $\WH_G$ (by Theorem \ref{convexhull}) and the orbit map from $G$ to $\WH_G$ is a quasi-isometry, this will complete the proof.

Let $\mathbb R$ be the set of real numbers equipped with its usual metric. We electrify $\mathbb R$ as follows. Let $\mathcal J(r)$ be the collection of all closed intervals $\mathbb J$ of length $r \geq 1$. For each $\mathbb J$ in $\mathcal J(r)$, we create a new point $w_{\mathbb J}$ and adjoin an interval of length $\frac 1 2$ between each point in $\mathbb J$ and $w_{\mathbb J}$.  We equip this space with the induced path metric and call the result $\Rel(r)$. It is not difficult to see that the inclusion $\mathbb R \to \Rel(r)$ is an $(r,0)$--quasi-isometry. 

Let $\tau$ be a bi-infinite geodesic in $\WH_G$. Identifying $\tau$ with $\mathbb R$, we consider the projection $\pi_\tau \co \T(S) \to \mathcal P (\tau)$ as a map $\pi_\tau \co \T(S) \to \mathcal P (\mathbb R)$ and, using the Axiom of Choice, we replace this with a map
\[
\Pi_\tau \co \T(S) \to \mathbb R 
\]
by demanding that $\Pi_\tau(X)$ be some element of $\pi_\tau(X)$, noting that for any two points $X$ and $Y$ in $\T(S)$, 
\[
\d_{\mathbb R}(\Pi_\tau(X), \Pi_\tau(Y)) \leq \mathrm{diam}(\pi_\tau(X) \cup \pi_\tau(Y)).
\]

By Theorem 4.3 of \cite{minskycrelle}, there is a constant $B \geq 1$, depending only on $\epsilon_0$, such that 
\[
\mathrm{diam}(\Pi_\tau(\mathrm{thin}(\alpha, \epsilon_0))) \leq B.
\]
We extend the projection $\Pi_\tau$ to a projection 
\[
\widehat \Pi_\tau \co \Tel(S) \to \Rel(B)
\]
by demanding that the restriction $\widehat \Pi_\tau |_{\T(S)} = \Pi_\tau$, that for each $\alpha$ in $\C_0(S)$, $\widehat \Pi_\tau(u_\alpha)=w_{\mathbb J}$ for some interval $\mathbb J$ in $\mathcal J(B)$ containing $\Pi_\tau(\mathrm{thin}(\alpha, \epsilon_0))$---again using the Axiom of Choice---and that an electric edge between $u_\alpha$ and a point $X$ in $\mathrm{thin}(\alpha, \epsilon_0)$ be sent isometrically to an electric edge joining $w_{\mathbb J}$ and $\Pi_\tau(X)$.

We claim that this projection is coarsely Lipschitz. To see this, let $X$ and $Y$ be elements of $\Tel(S)$ with $\d_{\Tel}(X,Y) \leq 1$. An electric geodesic $\gamma$ joining $X$ and $Y$ is a concatenation of Teichm\"uller geodesic segments in $\T(S)$ and paths in unions of electric edges. Since $\d_{\Tel}(X,Y) \leq 1$, we may write $\gamma$ as a concatenation of paths \linebreak $[X,X']\cup[X',Y']\cup[Y',Y]$, where $[X,X']$ and $[Y',Y]$ are paths in unions of electric edges, and $[X',Y']$ is a Teichm\"uller geodesic segment---any of which may be a constant path. 

Since $\WH_G$ is $\epsilon_0$--cobounded, so is $\tau$.  By Theorem \ref{quasiproj}, for any two points $X$ and $Y$ in $\T(S)$, 
\[
\d_{\mathbb R}(\Pi_\tau(X) , \Pi_\tau(Y)) \leq \d_\T(X,Y) + 4b,
\]
where $b =b(\epsilon_0)$, and since the inclusion $\mathbb R \to \Rel(B)$ is $1$--Lipschitz, we have
\[
\d_{\Rel}(\widehat \Pi_\tau(X'), \widehat \Pi_\tau(Y')) \leq 1 + 4b.
\]
Since $\widehat \Pi_\tau$ conducts electricity, we have
\[
\d_{\Rel}(\widehat \Pi_\tau(Z), \widehat \Pi_\tau(Z')) \leq 1
\]
when $Z \in \{X,Y\}$.
So, 
\[
\d_{\Rel}(\widehat \Pi_\tau(X), \widehat \Pi_\tau(Y)) \leq 3+4b
\]
and we see that $\widehat \Pi_\tau$ is $(C',C')$--Lipschitz, for $C' = 3 + 4b$.

Since $\mathbb R \to \Rel(B)$ is a $(B,0)$--quasi-isometry, for all $X$ and $Y$ in $\tau$
\begin{align*}
\d_{\Tel}(X,Y) & \geq \frac{1}{C'} \, \d_{\Rel}(X,Y) - 1 \\
& \geq \frac{1}{BC'} \, \d_{\mathbb R}(X,Y) - 1 \\
& = \frac{1}{BC'} \, \d_{\T}(X,Y) - 1.
\end{align*}
Letting $K = BC'$ and $C = 1$ completes the proof.
\end{proof}

\section{Questions}\label{questions}

\subsection{The analogy}\label{analogy}

For a Kleinian group, acting cocompactly on the domain of discontinuity is insufficient to guarantee convex cocompactness.  For example, L. Bers first established the existence of \textbf{singly degenerate} Kleinian groups isomorphic to the fundamental group of a hyperbolic surface \cite{bers}: geometrically infinite groups whose domains of discontinuity are topological disks on which the groups act cocompactly.

When drawing an analogy between $\T(S)$ and $\mathbb H^3$ it is in many respects prudent to compare $\Mod(S)$ with a Kleinian group $\Gamma$ of finite covolume. In this picture, the moduli space $\M(S)$ plays the role of the orbifold $M_\Gamma = \mathbb H^3 / \Gamma$ and as $\M(S)$ is non-compact---and has finite volume in a certain sense (compare Masur \cite{masurinterval}, Section 5)---the analogy suggests that $M_\Gamma$ also be non-compact.

The resolution of Marden's Tameness Conjecture by I. Agol \cite{agol} and (independently) D. Calegari and D. Gabai \cite{calgabai}, combined with R. Canary's Covering Theorem \cite{canary}, implies that a finitely generated subgroup of $\Gamma$ is either geometrically finite or virtually the fiber subgroup of a hyperbolic $3$--manifold fibering over the circle---see \cite{canarysurvey}.  In particular, no groups like the ones constructed by Bers can occur in $\Gamma$.

A cocompact action of a subgroup of $\Gamma$ on its domain of discontinuity is still however insufficient to guarantee convex cocompactness.
One can construct examples of subgroups $\Gamma_0 < \Gamma$ which are geometrically finite, but for which all cusps are rank 2 (and hence ``internal'' to the convex core).
In this situation, the convex hull of the limit set of $\Gamma_0$ is not cobounded with respect to the covering $\mathbb H^3 \to M_\Gamma$.
As we have seen, this behavior does not present itself in $\Mod(S)$, in light of Theorem \ref{domainsub}.  So, for subgroups of $\Mod(S)$ it is feasible that a cocompact action on the domain of discontinuity is equivalent to convex cocompactness.

\begin{question} If a finitely generated subgroup $G$ of $\Mod(S)$ acts cocompactly on $\Delta_G \neq \emptyset$, is it convex cocompact?
\end{question}
\noindent By Theorem \ref{domainsub}, an affirmative answer to this question would follow from an affirmative answer to the following

\begin{question}\label{hullquestion} If $G$ is a finitely generated subgroup of $\Mod(S)$ and $\WH_G$ is cobounded, is $G$ convex cocompact?
\end{question}

When $G$ is convex cocompact, $\Delta_G$ is the largest open set in $\PML(S)$ on which $G$ acts properly discontinuously.  This is also true for Veech groups.  

\begin{question} Let $G$ be a finitely generated subgroup of $\Mod(S)$.  Is $\Delta_G$ the largest open set in $\PML(S)$ on which $G$ acts properly discontinuously?
\end{question}

We note that the action of $G$ on the preimage of $\Delta_G$ in $\ML(S)$ is also properly discontinuous, and it has been shown by C. Lecuire \cite{lecuire} that the handlebody group $G$ (in genus at least 3) acts properly discontinuously on a strictly larger set in $\ML(S)$ than the preimage of $\Delta_G$.

\begin{question} If one takes $\check{\Delta}_G$ to be the largest open set on which $G$ acts properly discontinuously, is it true that $G$ is convex cocompact if and only if $(\T(S) \cup \check{\Delta}_G)/G$ is compact?
\end{question}

We note that the answer to this question is affirmative if the answer to the previous one is as well.

\subsection{Examples}\label{examples}

At present, the only known examples of convex cocompact subgroups of $\Mod(S)$ are virtually free.  To the authors' knowledge, the only known examples are: groups obtained by taking powers of independent pseudo-Anosov mapping classes; certain free products of finite subgroups of $\Mod(S)$, constructed by Honglin Min \cite{min}; 
and purely pseudo-Anosov subgroups of graphs of Veech groups, due to the second author \cite{leininger}.

In \cite{masurhandle}, Masur studies the group of mapping classes of $S$ that extend over a handlebody, called the \textbf{handlebody group}.

\begin{question} Is every finitely generated purely pseudo-Anosov subgroup of the handlebody group convex cocompact?
\end{question}

Let $\dot S$ denote the surface $S$ minus a point.  There is a short exact sequence
\[
1 \to \pi_1(S) \to \Mod(\dot S) \to \Mod(S) \to 1
\]
where an element of $\pi_1(S)$ is sent to the mapping class that ``spins'' the puncture about the corresponding loop in $S$ and $\Mod(\dot S) \to \Mod(S)$ is the map forgetting the puncture---see \cite{birman}.

\begin{question} Is every finitely generated purely pseudo-Anosov subgroup of $\pi_1(S)$ a convex cocompact subgroup of $\Mod(\dot S)$?
\end{question}

An affirmative answer to this question would show that K. Whittlesey's group \cite{whittle} is locally convex cocompact---this is a normal purely pseudo-Anosov subgroup of the mapping class group of a surface of genus two and is isomorphic to a free group of infinite rank.

\begin{question}[Farb--Mosher \cite{FMcc}] Is every finitely generated subgroup of Whittlesey's group convex cocompact?
\end{question}

A more delicate question is 

\begin{question}[Farb--Mosher \cite{FMcc}] Is there a convex cocompact subgroup $G$ of $\Mod(S)$ that is not virtually free?
\end{question}

And a more daring question is

\begin{question}[Reid \cite{reid}] Let $m \geq 3$ be less than the virtual cohomological dimension of $\Mod(S)$ and let $\Gamma$ be a torsion free uniform lattice in $\mathrm{SO}(m,1)$.  Is there an injection $\Gamma \hookrightarrow \Mod(S)$ whose image is purely pseudo-Anosov?
\end{question}

Note that if $\Gamma$ is the fundamental group of a closed fibered hyperbolic $3$--manifold with fiber subgroup $\Sigma$ and $\Gamma$ injects into $\Mod(S)$ with convex cocompact image, then $\Sigma$ could not act cocompactly on its weak hull $\WH_\Sigma$, as $\WH_\Sigma$ would equal $\WH_\Gamma$. Such a $\Sigma$ would resolve Question \ref{hullquestion} in the negative and it follows from work in \cite{FMcc} that the associated $\pi_1(S)$--extension of $\Sigma$ would be a non-hyperbolic group with a finite Eilenberg--\linebreak Mac Lane space and no Baumslag--Solitar subgroups---see \cite{klsurvey} and Question 1.1 of \cite{bestvinaproblems}.

\subsection{The sociology of ending laminations}\label{sociology}

Theorem \ref{complextheorem} implies that the Gromov boundary of a convex cocompact $G$ embeds in the boundary of $\C(S)$, the space $\EL(S)$ of potential ending laminations for hyperbolic $3$--manifolds homeomorphic to $S \times \mathbb R$.  So, if $\EL(S)$ is totally disconnected, then every convex cocompact subgroup of $\Mod(S)$ is virtually free.  To see this, note that $\partial G$ is compact and so, provided $G$ is not virtually cyclic, total disconnectedness of $\EL(S)$ along with the above embedding implies that $\partial G$ is a Cantor set.  Such a group is virtually free \cite{stallingsends,gromov,apresgromov}.

With this in mind, we close with a question of Peter Storm and a related question.

\begin{question}[Storm] Is $\EL(S)$ connected?  Is it path connected?
\end{question}

This is closely related to connectivity outside large balls in $\C(S)$.  Specifically, the following is unknown in general.

\begin{question} Does there exists an $A>0$ such that given any $R>0$ and any two points $u,v \in \C(S)$ outside a ball of radius $R$, there is a path connecting $u$ to $v$ that remains outside the ball (with the same center) of radius $R-A$?
\end{question}

The answer to this question has been resolved by S. Schleimer \cite{salty} when $S$ is a once-punctured surface with genus at least $2$.  Indeed, Schleimer shows in this case that the complement of any $R$-ball is path connected.

\bibliographystyle{plain}
\bibliography{convex}

\def\cprime{$'$}
\begin{thebibliography}{10}

\bibitem{abikoff}
William Abikoff.
\newblock {\em The real analytic theory of {T}eichm\"uller space}, volume 820
  of {\em Lecture Notes in Mathematics}.
\newblock Springer, Berlin, 1980.

\bibitem{agol}
Ian Agol.
\newblock {Tameness of hyperbolic 3-manifolds, Preprint}.
\newblock \texttt{arXiv:math.GT/0405568},.

\bibitem{ahlforsqc}
Lars~V. Ahlfors.
\newblock {\em Lectures on quasiconformal mappings}.
\newblock Manuscript prepared with the assistance of Clifford J. Earle, Jr. Van
  Nostrand Mathematical Studies, No. 10. D. Van Nostrand Co., Inc., Toronto,
  Ont.-New York-London, 1966.

\bibitem{bers}
Lipman Bers.
\newblock On boundaries of {T}eichm\"uller spaces and on {K}leinian groups.
  {I}.
\newblock {\em Ann. of Math. (2)}, 91:570--600, 1970.

\bibitem{bestvinaproblems}
Mladen Bestvina.
\newblock {Questions in geometric group theory}.
\newblock \texttt{http://www.math.utah.edu/$\sim$bestvina}.

\bibitem{combinationtheorem}
Mladen Bestvina and Mark Feighn.
\newblock A combination theorem for negatively curved groups.
\newblock {\em J. Differential Geom.}, 35(1):85--101, 1992.

\bibitem{birman}
Joan~S. Birman.
\newblock Mapping class groups and their relationship to braid groups.
\newblock {\em Comm. Pure Appl. Math.}, 22:213--238, 1969.

\bibitem{bonahon}
Francis Bonahon.
\newblock The geometry of {T}eichm\"uller space via geodesic currents.
\newblock {\em Invent. Math.}, 92(1):139--162, 1988.

\bibitem{bonahonsurvey}
Francis Bonahon.
\newblock Geodesic laminations on surfaces.
\newblock In {\em Laminations and foliations in dynamics, geometry and topology
  (Stony Brook, NY, 1998)}, volume 269 of {\em Contemp. Math.}, pages 1--37.
  Amer. Math. Soc., Providence, RI, 2001.

\bibitem{bowditchgeom}
Brian~H. Bowditch.
\newblock Geometrical finiteness for hyperbolic groups.
\newblock {\em J. Funct. Anal.}, 113(2):245--317, 1993.

\bibitem{bowditch}
Brian~H. Bowditch.
\newblock Intersection numbers and the hyperbolicity of the curve complex.
\newblock {\em J. Reine Angew. Math.}, 598:105--129, 2006.

\bibitem{BH}
Martin~R. Bridson and Andr{\'e} Haefliger.
\newblock {\em Metric spaces of non-positive curvature}, volume 319 of {\em
  Grundlehren der Mathematischen Wissenschaften [Fundamental Principles of
  Mathematical Sciences]}.
\newblock Springer-Verlag, Berlin, 1999.

\bibitem{brockfarb}
Jeffrey Brock and Benson Farb.
\newblock {Curvature and rank of Teichm{\"u}ller space}.
\newblock Preprint, \texttt{arXiv:math.GT/0109045}.

\bibitem{calgabai}
Danny Calegari and David Gabai.
\newblock {Shrinkwrapping and the taming of hyperbolic 3-manifolds}.
\newblock Preprint, \texttt{arXiv:math.GT/0407161}.

\bibitem{notesonnotes}
R.~D. Canary, D.~B.~A. Epstein, and P.~Green.
\newblock Notes on notes of {T}hurston.
\newblock In {\em Analytical and geometric aspects of hyperbolic space
  (Coventry/Durham, 1984)}, volume 111 of {\em London Math. Soc. Lecture Note
  Ser.}, pages 3--92. Cambridge Univ. Press, Cambridge, 1987.

\bibitem{canarysurvey}
Richard~D. Canary.
\newblock Covering theorems for hyperbolic {$3$}-manifolds.
\newblock In {\em Low-dimensional topology (Knoxville, TN, 1992)}, Conf. Proc.
  Lecture Notes Geom. Topology, III, pages 21--30. Internat. Press, Cambridge,
  MA, 1994.

\bibitem{canary}
Richard~D. Canary.
\newblock A covering theorem for hyperbolic {$3$}-manifolds and its
  applications.
\newblock {\em Topology}, 35(3):751--778, 1996.

\bibitem{coor}
M.~Coornaert, T.~Delzant, and A.~Papadopoulos.
\newblock {\em G\'eom\'etrie et th\'eorie des groupes}, volume 1441 of {\em
  Lecture Notes in Mathematics}.
\newblock Springer-Verlag, Berlin, 1990.

\bibitem{duchin}
Moon Duchin.
\newblock Thin triangles and a multiplicative ergodic theorem for
  {T}eichm{\"u}ller geometry.
\newblock Preprint, \texttt{arXiv:math.GT/0508046}.

\bibitem{FMcc}
Benson Farb and Lee Mosher.
\newblock Convex cocompact subgroups of mapping class groups.
\newblock {\em Geom. Topol.}, 6:91--152 (electronic), 2002.

\bibitem{FMII}
Benson Farb and Lee Mosher.
\newblock The geometry of surface-by-free groups.
\newblock {\em Geom. Funct. Anal.}, 12(5):915--963, 2002.

\bibitem{FLP}
A.~Fathi, F.~Laudenbach, and V.~Po\'enaru.
\newblock {\em Travaux de {T}hurston sur les surfaces}.
\newblock Soci\'et\'e Math\'ematique de France, Paris, 1991.
\newblock S\'eminaire Orsay, Reprint of {\it Travaux de Thurston sur les
  surfaces}, Soc.\ Math.\ France, Paris, 1979 Ast\'erisque No. 66-67 (1991).

\bibitem{gardiner}
Frederick~P. Gardiner.
\newblock {\em Teichm\"uller theory and quadratic differentials}.
\newblock Pure and Applied Mathematics (New York). John Wiley \& Sons Inc., New
  York, 1987.
\newblock A Wiley-Interscience Publication.

\bibitem{apresgromov}
{\'E}tienne Ghys and Pierre de~la Harpe, editors.
\newblock {\em Sur les groupes hyperboliques d'apr\`es {M}ikhael {G}romov},
  volume~83 of {\em Progress in Mathematics}.
\newblock Birkh\"auser Boston Inc., Boston, MA, 1990.
\newblock Papers from the Swiss Seminar on Hyperbolic Groups held in Bern,
  1988.

\bibitem{gromov}
Mikha{\"{i}}l Gromov.
\newblock Hyperbolic groups.
\newblock In {\em Essays in group theory}, volume~8 of {\em Math. Sci. Res.
  Inst. Publ.}, pages 75--263. Springer, New York, 1987.

\bibitem{halpern}
Noemi Halpern.
\newblock A proof of the collar lemma.
\newblock {\em Bull. London Math. Soc.}, 13(2):141--144, 1981.

\bibitem{hamenstadt}
Ursula Hamenst{\"a}dt.
\newblock {Word hyperbolic extensions of surface groups}.
\newblock Preprint, \texttt{arXiv:math.GT/0505244}.

\bibitem{hubbardmasur}
John Hubbard and Howard Masur.
\newblock Quadratic differentials and foliations.
\newblock {\em Acta Math.}, 142(3-4):221--274, 1979.

\bibitem{klsurvey}
Richard P.~Kent IV and Christopher~J. Leininger.
\newblock {Subgroups of the mapping class group from the geometrical viewpoint
  }.
\newblock To appear in \textit{Proceedings of the 2005 Ahlfors--Bers
  Colloquium}, \texttt{arXiv:math.GT/0702034}.

\bibitem{ivanov}
Nikolai~V. Ivanov.
\newblock {\em Subgroups of {T}eichm\"uller modular groups}, volume 115 of {\em
  Translations of Mathematical Monographs}.
\newblock American Mathematical Society, Providence, RI, 1992.
\newblock Translated from the Russian by E. J. F. Primrose and revised by the
  author.

\bibitem{keen}
Linda Keen.
\newblock Collars on {R}iemann surfaces.
\newblock In {\em Discontinuous groups and Riemann surfaces (Proc. Conf., Univ.
  Maryland, College Park, Md., 1973)}, pages 263--268. Ann. of Math. Studies,
  No. 79. Princeton Univ. Press, Princeton, N.J., 1974.

\bibitem{kerckhoff}
Steven~P. Kerckhoff.
\newblock The asymptotic geometry of {T}eichm\"uller space.
\newblock {\em Topology}, 19(1):23--41, 1980.

\bibitem{klarreich}
Erica Klarreich.
\newblock The boundary at infinity of the curve complex and the relative
  {T}eichm\"uller space.
\newblock Preprint, \texttt{http://nasw.org/users/klarreich/publications.htm},.

\bibitem{lecuire}
Cyril Lecuire.
\newblock {Structure hyperboliques convexes sur les vari{\'{e}}t{\'{e}}s de
  dimension 3}, \textit{th{\`{e}}se de doctorat de l'{E}{N}{S} {L}yon}.
\newblock \texttt{http://www.maths.warwick.ac.uk/$\sim$clecuire}.

\bibitem{leininger}
Christopher~J. Leininger.
\newblock Graphs of veech groups.
\newblock Work in progress.

\bibitem{levitt}
Gilbert Levitt.
\newblock Foliations and laminations on hyperbolic surfaces.
\newblock {\em Topology}, 22(2):119--135, 1983.

\bibitem{Mar}
Albert Marden.
\newblock The geometry of finitely generated kleinian groups.
\newblock {\em Ann. of Math. (2)}, 99:383--462, 1974.

\bibitem{masurclass}
Howard Masur.
\newblock On a class of geodesics in {T}eichm\"uller space.
\newblock {\em Ann. of Math. (2)}, 102(2):205--221, 1975.

\bibitem{masuruniquely}
Howard Masur.
\newblock Uniquely ergodic quadratic differentials.
\newblock {\em Comment. Math. Helv.}, 55(2):255--266, 1980.

\bibitem{masurinterval}
Howard Masur.
\newblock Interval exchange transformations and measured foliations.
\newblock {\em Ann. of Math. (2)}, 115(1):169--200, 1982.

\bibitem{twoboundaries}
Howard Masur.
\newblock Two boundaries of {T}eichm\"uller space.
\newblock {\em Duke Math. J.}, 49(1):183--190, 1982.

\bibitem{masurhandle}
Howard Masur.
\newblock Measured foliations and handlebodies.
\newblock {\em Ergodic Theory Dynam. Systems}, 6(1):99--116, 1986.

\bibitem{masurhaus}
Howard Masur.
\newblock Hausdorff dimension of the set of nonergodic foliations of a
  quadratic differential.
\newblock {\em Duke Math. J.}, 66(3):387--442, 1992.

\bibitem{MM1}
Howard~A. Masur and Yair~N. Minsky.
\newblock Geometry of the complex of curves. {I}. {H}yperbolicity.
\newblock {\em Invent. Math.}, 138(1):103--149, 1999.

\bibitem{MM2}
Howard~A. Masur and Yair~N. Minsky.
\newblock Geometry of the complex of curves. {II}. {H}ierarchical structure.
\newblock {\em Geom. Funct. Anal.}, 10(4):902--974, 2000.

\bibitem{MMunstable}
Howard~A. Masur and Yair~N. Minsky.
\newblock Unstable quasi-geodesics in {T}eichm\"uller space.
\newblock In {\em In the tradition of Ahlfors and Bers (Stony Brook, NY,
  1998)}, volume 256 of {\em Contemp. Math.}, pages 239--241. Amer. Math. Soc.,
  Providence, RI, 2000.

\bibitem{masurwolf}
Howard~A. Masur and Michael Wolf.
\newblock Teichm\"uller space is not {G}romov hyperbolic.
\newblock {\em Ann. Acad. Sci. Fenn. Ser. A I Math.}, 20(2):259--267, 1995.

\bibitem{mccarthypapa}
John McCarthy and Athanase Papadopoulos.
\newblock Dynamics on {T}hurston's sphere of projective measured foliations.
\newblock {\em Comment. Math. Helv.}, 64(1):133--166, 1989.

\bibitem{min}
Honglin Min.
\newblock Work in progress.

\bibitem{ELCI}
Yair~N. Minsky.
\newblock {The classification of Kleinian surface groups, I: Models and
  bounds}.
\newblock to appear \textit{{A}nn. of {M}ath.}, \texttt{arXiv:math.GT/0302208}.

\bibitem{minskycrelle}
Yair~N. Minsky.
\newblock Quasi-projections in {T}eichm\"uller space.
\newblock {\em J. Reine Angew. Math.}, 473:121--136, 1996.

\bibitem{minskyinj}
Yair~N. Minsky.
\newblock Kleinian groups and the complex of curves.
\newblock {\em Geom. Topol.}, 4:117--148, 2000.

\bibitem{minskybound}
Yair~N. Minsky.
\newblock Bounded geometry for {K}leinian groups.
\newblock {\em Invent. Math.}, 146(1):143--192, 2001.

\bibitem{Mitra1}
Mahan Mitra.
\newblock Ending laminations for hyperbolic group extensions.
\newblock {\em Geom. Funct. Anal.}, 7(2):379--402, 1997.

\bibitem{Mitra2}
Mahan Mitra.
\newblock Cannon-{T}hurston maps for hyperbolic group extensions.
\newblock {\em Topology}, 37(3):527--538, 1998.

\bibitem{Mitra3}
Mahan Mitra.
\newblock Cannon-{T}hurston maps for trees of hyperbolic metric spaces.
\newblock {\em J. Differential Geom.}, 48(1):135--164, 1998.

\bibitem{hypbyhyp}
Lee Mosher.
\newblock A hyperbolic-by-hyperbolic hyperbolic group.
\newblock {\em Proc. Amer. Math. Soc.}, 125(12):3447--3455, 1997.

\bibitem{mumford}
David Mumford.
\newblock A remark on {M}ahler's compactness theorem.
\newblock {\em Proc. Amer. Math. Soc.}, 28:289--294, 1971.

\bibitem{rafi}
Kasra Rafi.
\newblock A characterization of short curves of a {T}eichm\"uller geodesic.
\newblock {\em Geom. Topol.}, 9:179--202, 2005.

\bibitem{reid}
Alan~W. Reid.
\newblock Surface subgroups of mapping class groups.
\newblock In {\em Problems on mapping class groups and related topics},
  volume~74 of {\em Proc. Sympos. Pure Math.}, pages 257--268. Amer. Math.
  Soc., Providence, RI, 2006.

\bibitem{salty}
Saul Schleimer.
\newblock The end of the curve complex.
\newblock {Preprint, \texttt{arXiv:math.GT/0608505}},.

\bibitem{serretf}
J.-P. Serre.
\newblock Rigidit\'e de foncteur d'{J}acobi d'{\'e}chelon {$n \geq 3$}.
\newblock In {\em S{\'e}minaire Henri Cartan, 1960/61, Expos{\'e} 17
  (Appendice)}. Secr{\'e}tariat math{\'e}matique, Paris, 1960/61.

\bibitem{stallingsends}
John~R. Stallings.
\newblock On torsion-free groups with infinitely many ends.
\newblock {\em Ann. of Math. (2)}, 88:312--334, 1968.

\bibitem{Sul}
Dennis Sullivan.
\newblock Quasiconformal homeomorphisms and dynamics. {II}. {S}tructural
  stability implies hyperbolicity for {K}leinian groups.
\newblock {\em Acta Math.}, 155(3-4):243--260, 1985.

\bibitem{swenson}
Eric~L. Swenson.
\newblock Quasi-convex groups of isometries of negatively curved spaces.
\newblock {\em Topology Appl.}, 110(1):119--129, 2001.
\newblock Geometric topology and geometric group theory (Milwaukee, WI, 1997).

\bibitem{whittle}
Kim Whittlesey.
\newblock Normal all pseudo-{A}nosov subgroups of mapping class groups.
\newblock {\em Geom. Topol.}, 4:293--307 (electronic), 2000.

\end{thebibliography}

\bigskip

\noindent Department of Mathematics, Brown University, Providence, RI 02912 \newline \noindent \texttt{rkent@math.brown.edu}

\bigskip

\noindent Department of Mathematics, University of Illinois, Urbana-Champaign, IL 61801 \newline \noindent  \texttt{clein@math.uiuc.edu}

\end{document}